\documentclass[11pt]{amsart}
\usepackage{amssymb,amsmath,txfonts,mathrsfs,mathtools,titletoc}
\usepackage[alphabetic]{amsrefs}

\newtheorem{theorem}{Theorem}[section]
\newtheorem{prop}[theorem]{Proposition}
\newtheorem{lemma}[theorem]{Lemma}
\newtheorem{remark}[theorem]{Remark}
\newtheorem{claim}[theorem]{Claim}

\newtheorem{question}[theorem]{Question}

\newtheorem{definition}[theorem]{Definition} 
\newtheorem{cor}[theorem]{Corollary}
\newtheorem{example}[theorem]{Example}
\newtheorem{conj}[theorem]{Conjecture}

\numberwithin{equation}{section}

\def\pf{{\it Proof:}~}

\begin{document}

\title[Three circles theorems for harmonic functions]{Three circles theorems for harmonic functions}
\author{Guoyi Xu}
\address{Yau Mathematical Sciences Center\\Jin Chun Yuan West Building \\ Tsinghua University, Beijing\\P. R. China, 100084}
\email{gyxu@math.tsinghua.edu.cn}

\begin{abstract}
We proved two Three Circles Theorems for harmonic functions on manifolds in integral sense. As one application, on manifold with nonnegative Ricci curvature, whose tangent cone at infinity is the unique metric cone with unique conic measure,  we showed the existence of nonconstant harmonic functions with polynomial growth. This existence result recovered and generalized the former result of Y. Ding, and led to a complete answer of L. Ni's conjecture. Furthermore in similar context, combining the techniques of estimating the frequency of harmonic functions with polynomial growth, which were developed by Colding and Minicozzi, we confirmed their conjecture about the uniform bound of frequency. 
\end{abstract}


\maketitle

\titlecontents{section}[0em]{}{\hspace{.5em}}{}{\titlerule*[1pc]{.}\contentspage}
\titlecontents{subsection}[1.5em]{}{\hspace{.5em}}{}{\titlerule*[1pc]{.}\contentspage}
\tableofcontents


\section{Introduction}\label{SECTION 1}
In $1975$, S.-T. Yau \cite{Yau} generalized the classical Liouville theorem to complete manifolds with nonnegative Ricci curvature. Specially, he proved that any positive harmonic function on such manifolds is constant. In \cite{Cheng-sub}, S.-Y. Cheng further proved that on such manifolds any harmonic function of sublinear growth must be constant. On complete manifolds, harmonic functions with polynomial growth are important analytic functions besides the Green's function and the heat kernel (note the latter two have singularities). 

In the study of harmonic functions on complete manifolds, Yau considered the space of harmonic functions with polynomial growth:
\begin{align}
\mathscr{H}_d(M)=\{u|\ u(x) \ is\ harmonic\ on \ M^n\ ; |u(x)|\leq K(\rho(x)+ 1)^d\ for\ some \ K> 0\} \nonumber
\end{align}
where $\rho(x)= d(p, x)$ and $p$ is some fixed point on $M^n$. In \cite{problem}, the following conjecture was made:

\begin{conj}[S.-T. Yau]\label{conj 1.1}
{Let $M^n$ be a complete manifold with nonnegative Ricci curvature, then $dim\big(\mathscr{H}_d(M)\big)< \infty$ for any $d> 0$.
}
\end{conj}

P. Li and L.-F. Tam firstly proved the conjecture for linear growth harmonic function ($d= 1$) in \cite{LT1}, and they further verified the $2$-dimensional case ($n= 2$) in \cite{LT2}. In $1997$, this conjecture was completely proved by Colding and Minicozzi \cite{CMAnn} (also see \cite{Li2}, \cite{CMcpam}, \cite{CMinv} and \cite{LW} for further developments). Although Yau's conjecture was confirmed, there are still several important questions about harmonic functions with polynomial growth remained open. 

It is well known that on any complete noncompact manifold, there always exist nonconstant harmonic functions (see \cite{GW}). A natural question is about the existence of nonconstant harmonic function with polynomial growth. Note any complete manifold with nonnegative Ricci curvature has at least linear volume growth (see \cite{Yau}). C. Sormani proved the following:
\begin{theorem}[\cite{Sormani}]\label{thm Sormani}
{Let $M^n$ be a complete manifold with nonnegative Ricci curvature and at most linear volume growth. If there exists a nonconstant harmonic function of polynomial growth, then the manifold splits isometrically, $M^n= \mathbb{R}\times N^{n-1}$.
}
\end{theorem}
As observed in \cite{Ding2}, there exists a suitable metric on $\mathbb{R}^{+}\times \mathbb{S}^{n- 1}$, which has nonnegative Ricci curvature and linear volume growth, but can not split isometrically. From Theorem \ref{thm Sormani}, it will not admit any nonconstant harmonic functions with polynomial growth. A concrete example is given in Section \ref{SECTION 4} (Example \ref{example Ding's}). 

Hence, to study the existence of nonconstant harmonic function with polynomial growth in nonnegative Ricci curvature context, we need to have some restriction on the volume growth of $M^n$. To the author's knowledge, the following question is still open :
\begin{question}\label{ques 1.3}
{If $(M^n, g)$ is a complete manifold with $Rc\geq 0$ and maximal volume growth, does there exist $d\geq 1$ such that $dim\big(\mathscr{H}_d(M)\big)\geq 2$? In other words, is there any nonconstant harmonic function with polynomial growth order at most $d$ on $M^n$?
}
\end{question}

In another direction, putting the existence problem in the positively curved context, L. Ni \cite{Ni} made the following conjecture:
\begin{conj}[L. Ni]\label{conj 5.4}
{Let $(M^n, g)$ be a complete Riemannian manifold of positive sectional curvature, the necessary and sufficient condition that $M^n$ admits nonconstant harmonic functions of polynomial growth is, that $M^n$ is of maximum volume growth.
}
\end{conj}

For the corresponding conjecture on K\"ahler manifolds, L. Ni proved that if the manifold is complete K\"ahler with bounded nonnegative bisectional curvature and of maximum volume growth, it admits nonconstant holomorphic functions of polynomial growht (see Corollary $1$ of \cite{Ni-mrl}). On the other hand, recently, G. Liu (see Theorem $2$ of \cite{Liu2}) showed that if $M$ is a complete noncompact K\"ahler manifold with positive bisectional curvature, and it admits a nonconstant holomorphic function with polynomial growth, then $M$ is of maximal volume growth. 

In Riemannian geometry context, Y. Ding \cite{Ding 2} proved that on complete manifolds with $Rc\geq 0$, maximal volume growth and the unique tangent cone at infinity, there exists nonconstant harmonic function with polynomial growth. Note when $M^n$ has nonnegative sectional curvature, from Theorem $I. 26$ in \cite{RFTA}, the tangent cone at infinity of $M^n$ is the unique metric cone $C(Y)$, where $Y$ is a compact metric space. Hence Ding's existence result will imply the sufficient part of Conjecture \ref{conj 5.4}. However we have some difficulties to verify the proof of Ding's existence result. For example, Lemma $1.1$, Corollary $1.11$ and Lemma $1.2$ in \cite{Ding 2} does not hold for zero function.

The main technical tool in \cite{Ding 2} is a generalization of the monotonicity of frequency for harmonic functions on $\mathbb{R}^n$, which is a type of Three Circles Theorem in $L^2$ sense on complete manifolds (also see \cite{Zhang} and \cite{CDM} for related results). 

The classical Hadamard's Three Circles Theorem was implied in the announcement \cite{Hadamard} published in $1896$, we state it as the following form, which is sort of consistent with our presentation in Theorem \ref{thm 3.6 induction of J} and Theorem \ref{thm 3.5 induction of I} of this paper.
\begin{theorem}[Hadamard's Three Circles Theorem]\label{thm TCT}
{If $f(z)$ is a holomorphic function on $\{z|\ \frac{r}{4}\leq |z|\leq r\}\subset \mathbb{C}$, where $r> 0$ is some constant, then
\begin{align}
\frac{M(r/2)}{M(r/4)}\leq \frac{M(r)}{M(r/2)}  \label{equ TCT}
\end{align}
where $M(s)\vcentcolon =\max_{|z|= s}\{|f(z)|\}$, 
}
\end{theorem}

The classical Hadamard Three Circles Theorem for holomorphic functions had also been generalized to solutions of partial differential equations in different contexts by L. Simon \cite{Simon}, J. Cheeger and G. Tian \cite{CT}, G. Liu \cite{Liu}. In spite of our concerns about the argument in \cite{Ding2}, partially motivated by the results there, we proved two modified Three Circles Theorems for harmonic functions in integral sense (see Theorem \ref{thm 3.6 induction of J} and Theorem \ref{thm 3.5 induction of I} in Section \ref{SECTION 3}).  

\begin{remark}\label{rem three circles}
{One novel thing in one of our Three Circles theorems (Theorem \ref{thm 3.6 induction of J}) is that it can be applied for the collapsed case, i.e. the case that the maximal volume growth assumption does not hold. Also Theorem \ref{thm 3.6 induction of J} is dealing with the integral $J_u$, whose domain is different from the original one considered in \cite{Ding 2} (see Corollary $1.11$ there). View from the domains of the integral appearing in those results, Theorem \ref{thm 3.6 induction of J} is more like a three disks theorem, and Ding's technical tool is based on a three annuli theorem.
}
\end{remark}

As one application of our Three Circles Theorems, we proved the following theorem, which generalizes the existence result in \cite{Ding2}.

\begin{theorem}\label{thm existence of poly growth harmonic function}
{Let $(M^n, g)$ be a complete manifold with nonnegative Ricci curvature, the tangent cone at infinity with renormalized limit measure is a unique metric cone $C(X)$ with the unique conic measure of power $\kappa\geq 2$, and $\mathcal{H}^{1}(X)> 0$. Then 
\begin{align}
\inf \{\alpha|\ \alpha\in \mathscr{D}(M), \alpha\neq 0\}< \infty \label{upper bound of ds}
\end{align}
and for any $d> \inf \{\alpha|\ \alpha\in \mathscr{D}(M), \alpha\neq 0\}$,
\begin{align}
dim\big(\mathscr{H}_d(M)\big)\geq 2 \label{lower bound of dim}
\end{align}
}
\end{theorem}
The conic measure of power $\kappa$ and $\mathscr{D}(M)$ will be defined in Section \ref{SECTION 2}.

Generally, we do not know the uniqueness of the renormalized limit measure with respect to one tangent cone at infinity of the manifold (compare Example $1.24$ of \cite{CC1}). However, from Theorem $5.9$ of \cite{CC1}, if $(M^n, g)$ is a complete manifold with $Rc\geq 0$ and maximal volume growth, then every tangent cone at infinity has its unique renormalized measure; and the measure is a multiple of Hausdorff measure $\mathcal{H}^n$, which is a conic measure of power $n$. Hence, Theorem \ref{thm existence of poly growth harmonic function} implies the existence result in \cite{Ding2} mentioned above. The uniqueness of tangent cone at infinity is an important and hard problem, which was addressed in \cite{CT} and \cite{CMunique} for Ricci flat manifolds under various assumptions.

\begin{remark}\label{rem sharp}
{When $(M^n, g)$ is a complete manifold with $Rc\geq 0$ and maximal volume growth, S. Honda \cite{Honda} showed $dim\big(\mathscr{H}_d(M)\big)= 1$ for any $d< \inf \{\alpha|\ \alpha\in \mathscr{D}(M), \alpha\neq 0\}$. Hence the conclusion of Theorem \ref{thm existence of poly growth harmonic function} seems to be a sharp result except that the critical case $d= \inf \{\alpha|\ \alpha\in \mathscr{D}(M), \alpha\neq 0\}$ is not clear yet.
}
\end{remark}

To prove Theorem \ref{thm existence of poly growth harmonic function}, we partially followed the strategy in \cite{Ding2}. Because our Three Circles Theorem (Theorem \ref{thm 3.6 induction of J}) works for the collapsed case too, we succeeded in proving the existence of nonconstant harmonic functions with polynomial growth in the collapsed case first time, although with additional assumptions.

More concretely, to construct one nonconstant harmonic function with polynomial growth, we firstly choose a suitable harmonic function of polynomial growth on the tangent cone at infinity of the manifold, where we used the assumptions that the manifold has the conic renormalized limit measure and the tangent cone at infinity is a metric cone. 

Then we use the results of Cheeger \cite{Cheeger}, to get a sequence of approximate functions defined on a sequence of increasing geodesic balls exhausting the manifold, which are vanishing at the same fixed point. Solving the Dirichlet problem on those geodesic balls with the same boundary conditions as the corresponding approximate functions, this yields a family of harmonic functions defined on the exhausting domains of the manifold.

Because the sequence of harmonic functions constructed as above have the asymptotic growth behavior as the chosen harmonic function of polynomial growth on the tangent cone at infinity, we can get that the the ratio between the average integrals of those harmonic functions on bigger domain and smaller domain are uniformly bounded near infinity, where the bound depends on the growth rate of the harmonic function of polynomial growth on the tangent cone at infinity chosen above.

If we can get \textbf{`some induction estimate'} of the ratios from outer domains to inner domains, the uniform polynomial bound of the family of harmonic functions will be obtained by the induction method. Then, after the suitable rescaling, using the well-known Cheng-Yau's gradient estimate for harmonic functions in \cite{CY}, combining with the Arzela-Ascoli theorem, for some subsequence of these harmonic functions, we get the limit function defined on the whole manifold, which is harmonic function of polynomial growth. The nonconstancy of the limit function follows from its vanishing at the fixed point, and the non-vanishing of some local integral of the limit function, which resulted from the suitable chosen rescaling mentioned above.

Our Three Circles Theorem (Theorem \ref{thm 3.6 induction of J}) will play the role of the \textbf{`induction estimate'} needed in the above argument. Starting from the eigenfunctions expansion of harmonic functions on the metric cone with conic measure, the key idea to prove the Three Circles Theorem, is to use the gap between the eigenvalues of the tangent cone's cross-section. When the tangent cone at infinity with renormalized limit measure is the unique metric cone with unique conic measure, this gap is implied by the discreteness of the spectrum of Laplace operator on the cross-section, and we get the last piece in the proof of Theorem \ref{thm existence of poly growth harmonic function}.

When nonconstant harmonic functions with polynomial growth exist on complete manifolds, as proved in \cite{CM-JDG-97}, the bound of frequency is essential to describe the asymptotic structure of those functions (like the almost separation of variables). Hence a natural question is about the uniform bound of frequency of harmonic functions on manifolds. Based on the study in \cite{CM-JDG-97}, Colding and Minicozzi posed the following conjecture:
\begin{conj}[\cite{CM}]\label{conj 1.2}
{Suppose that $M^n$ has nonnegative Ricci curvature and maximal volume growth. If $u\in \mathscr{H}_d(M)$ for some $d> 0$, then the frequency of $u$ is uniformly bounded.
}
\end{conj}

\begin{remark}\label{rem history of fre}
{Besides \cite{CM-JDG-97} and \cite{CM}, the frequency was also studied in \cite{Alm} and \cite{GL}. For more related reference about the frequency, the reader can consult Remark $2.16$ in \cite{CM-JDG-97}. We would like to point out that the monotonicity of the frequency in Euclidean space can be viewed as the quantitative version of the classical Three Circles Theorem. 
}
\end{remark}

Roughly say, to get the uniform bound of the frequency, we only need to control the ratios of $I(r)$ on concentric circles with increasing radii. Checking the results and techniques developed in \cite{CM-JDG-97} carefully, the ratios have uniform bound on a sequence of concentric circles, whose radii are approaching the infinity. If a suitable Three Circles Theorem is available, the uniform bound of ratios can be obtained by the induction method similar as the former argument, which will imply the uniform bound of the frequency. Hence, using the Three Circles Theorem (Theorem \ref{thm 3.5 induction of I}) established in Section \ref{SECTION 3}, we proved the following theorem:

\begin{theorem}\label{thm unique cone case}
{Suppose that $(M^n, g)$ has nonnegative Ricci curvature and maximal volume growth, also assume the tangent cone at infinity of $M^n$ is unique. Then for $u(x)\in \mathscr{H}_d(M)$, the frequency of $u(x)$ is uniformly bounded by $C(u, n, V_M, d)$.
}
\end{theorem}

\begin{remark}
{This theorem confirms Conjecture \ref{conj 1.2} with the additional assumption the uniqueness of the tangent cone at infinity of manifolds. In fact, we proved a stronger result which implies Theorem \ref{thm unique cone case}, see Theorem \ref{thm 4.5 frequency of linear growth function} in Section \ref{SECTION 5} for details.
}
\end{remark}

As we mentioned before, from \cite{RFTA}, for any complete manifold with nonnegative sectional curvature, the tangent cone at infinity is a unique metric cone. From Theorem \ref{thm existence of poly growth harmonic function}, we have the following corollary.
\begin{cor}\label{cor disprove Ni's conjecture}
{Suppose that $(M^n, g)$ has nonnegative sectional curvature, the tangent cone at infinity with renormalized limit measure is a unique metric cone $C(X)$ with the unique conic measure of power $\kappa\geq 2$, and $\mathcal{H}^{1}(X)> 0$.  Then (\ref{upper bound of ds}) and (\ref{lower bound of dim}) hold.
}
\end{cor}

\begin{remark}\label{rem Ni's conjecture}
{On non-negatively curved manifolds, maximal volume growth implies the uniqueness of the tangent cone at infinity and the conic renormalized limit measure of power $\kappa= n$ and $\mathcal{H}^{n- 1}(X)> 0$, hence the sufficient part of Conjecture \ref{conj 5.4} is implied by Corollary \ref{cor disprove Ni's conjecture}. 

On the other hand, there exists complete manifolds $M^n$ with positive sectional curvature, whose tangent cone at infinity $C(X)$ has the unique renormalized limit measure, which is conic measure of power $\kappa\geq 2$ and $\mathcal{H}^{1}(X)> 0$ (see Example \ref{example Ni's}). Hence by Corollary \ref{cor disprove Ni's conjecture}, Example \ref{example Ni's} is a counterexample to the necessary part of Conjecture \ref{conj 5.4}. 
}
\end{remark}

The organization of this paper is as the following. In Section \ref{SECTION 2}, we stated some background facts about Gromov-Hausdorff convergence and Cheeger-Colding's theory, which are needed for later sections. We also recalled the definition of frequency function and the related formulas. 

In Section \ref{SECTION 3}, we proved two Three Circles Theorems, which are the key technical tools applicable for the existence and frequency problems respectively. For both theorems, the method is proof by contradiction and reduced the related analysis to the analysis on the tangent cone at infinity.

In Section \ref{SECTION 4}, we constructed the nonconstant harmonic function of polynomial growth from the harmonic function on the tangent cones at infinity. And the Three Circles Theorem is used to guarantee the polynomial growth of the constructed harmonic function. We also constructed two example manifolds, which address the nonexistence and existence of harmonic functions with polynomial growth, under linear volume growth and at least quadratic volume growth assumptions respectively. Specially, one example is the first counterexample to the necessary part of Conjecture \ref{conj 5.4}. 
 
In Section \ref{SECTION 5}, using the other Three Circles Theorem, combining the results and techniques developed in \cite{CM-JDG-97}, we proved the uniform bound of frequency. Some technical results in this section are well-known from \cite{CM-JDG-97} in more general context, but we provide the details here to make our argument self-contained in this concrete case.


\section{Background and notations}\label{SECTION 2}

In this section, we always assume that $(M^n, g)$ is an $n$-dimensional complete manifold with $Rc\geq 0$. We firstly review some background material about Gromov-Hausdorff convergence and analysis on limit spaces. 

Let $\big\{(M_i^n, p_i, \rho_i)\big\}$ be a sequence of pointed Riemannian manifolds, where $p_i\in M_i^n$ and $\rho_i$ is the metric on $M_i^n$. If $\big\{(M_i^n, p_i, \rho_i)\big\}$ converges to $(M_{\infty}, p_{\infty}, \rho_{\infty})$ in the Gromov-Hausdorff sense, we write $\displaystyle (M_i^n, p_i, \rho_i)\stackrel{d_{GH}}{\longrightarrow} (M_{\infty}, p_{\infty}, \rho_{\infty})$. See \cite{Gromov} for the definition and basic facts concerning Gromov-Hausdorff convergence.

A metric space $(M_{\infty}, p_{\infty}, \rho_{\infty})$ is a \textbf{tangent cone at infinity} of $M^n$ if it is a Gromov-Hausdorff limit of a sequence of rescaled manifolds $(M^n, p, r_j^{-2}g)$, where $r_j\rightarrow \infty$. By Gromov's compactness theorem, \cite{Gromov}, any sequence $r_j\rightarrow \infty$, has a subsequence, also denoted as $r_j\rightarrow \infty$, such that the rescaled manifolds $(M^n, p, r_j^{-2}g)$ converge to some tangent cone at infinity $M_{\infty}$ in the Gromov-Hausdorff sense. 

Let us recall that from Bishop-Gromov's volume comparison theorem, we can define the asymptotic volume ratio
\begin{align}
V_{M}= \lim_{r\rightarrow \infty}\frac{V(r)}{r^n} \label{AVR}
\end{align}
where $V(r)$ is the volume of the geodesic ball $B(r)$ centered at $p$ with radius $r$. And the above definition is independent of $p$, so we omit $p$ there. If $V_{M}> 0$, we say that $(M^n, g)$ has \textbf{maximal volume growth}. Note $V_M\leq V_0^n(1)$ from Bishop-Gromov's volume comparison theorem, where $V_{k}^{n}(r)$ is the volume of ball with radius $r$ in the $n$-dimensional space form with sectional curvature equal to $k$.

Example of Perelman \cite{Pere} shows that tangent cone at infinity is not unique in general even if the manifold with $Rc\geq 0$ has maximal volume growth and quadratic curvature decay. Although the tangent cone at infinity may be not unique, under maximal volume growth assumption, Cheeger and Colding proved the following theorem characterizing it:
\begin{theorem}[\cite{CC}]\label{thm metric cone at infinity}
{Let $M^n$ be a complete manifold with $Rc\geq 0$ and maximal volume growth, then every tangent cone at infinity $M_{\infty}$ is a metric cone $C(X)$, where $X$ is a compact metric space and $diam(X)\leq \pi$.
}
\end{theorem}
Note the metric on the metric cone $C(X)$ is $dr^2+ r^2dX$, where $r\in [0, \infty)$.

In the collapsed case (i.e. the maximal volume growth assumption does not hold), we can consider the renormalized measure on the limit space under the measured Gromov-Hausdorff convergence. As in Section $9$ of \cite{Cheeger}, we have the following definition.
\begin{definition}\label{def MGH}
{If $\omega_i$, $\omega_{\infty}$ are Borel regular measures on $M_i^n$, $M_{\infty}$, we say that $(M_i^n, p_i, \rho_i, \omega_i)$ converges to $(M_{\infty}, p_{\infty}, \rho_{\infty}, \omega_{\infty})$ in the \textbf{measured Gromov-Hausdorff sense}, if $(M_i^n, p_i, \rho_i)\stackrel{d_{GH}}{\longrightarrow} (M_{\infty}, p_{\infty}, \rho_{\infty})$, in addition, for any $x_i\rightarrow x_{\infty}$, ($x_i\in M_i^n$, $x_{\infty}\in M_{\infty}$), $r> 0$, we have 
\[\omega_{i}\Big(B_{i}(x_i, r)\Big)\rightarrow \omega_{\infty}\Big(B_{\infty} (x_{\infty}, r)\Big)\]
where $(M_{\infty}, \rho_{\infty})$ is a length space with length metric $\rho_{\infty}$, and  
\begin{align}
B_i(x_{i}, r)= \{z\in M_i^n|\  d_{\rho_i}(z, x_i)\leq r\}\ , \quad B_{\infty}(x_{\infty}, r) = \{z\in M_{\infty}|\  d_{\rho_{\infty}}(z, x_{\infty})\leq r\} \nonumber
\end{align}  
}
\end{definition}

For later use, we also set up the following \textbf{Blow Down Setup}: Note that $(M^n, g, \mu)$ is a complete Riemannian manifold with $Rc\geq 0$, where $\mu$ is the volume element determined by the metric $g$. We can define $(M_{i}, p, \rho_i,  \nu_i)$, where $M_i$ is the same differential manifold as $M^n$, $\rho_i$ is the metric defined as $\rho_i= r_i^{-2}g$, $\{r_i\}_{i= 1}^{\infty}$ is an increasing positive sequence whose limit is $\infty$, $p$ is a fixed point on $M_i= M^n$, and $\nu_i$ is the \textbf{renormalized measure} defined by 
\begin{align}
\nu_i (A)\vcentcolon= \Big(\int_{B_i(1)} 1 d\mu_{i}\Big)^{-1} \Big(\int_{A} 1 d\mu_{i}\Big)= r_i^{n} V(r_i)^{-1}\mu_i(A) \label{def of nu}
\end{align}
where $A\subset M_i$, $B_i(1)\vcentcolon= \{z\in M_i|\ d_{\rho_i}(z, p)\leq 1\}$, and $\mu_{i}$ is the volume element determined by $\rho_i$. Then by Gromov's compactness theorem (see \cite{Gromov}) and Theorem $1.6$ in \cite{CC1}, after passing to a suitable subsequence, we have
$(M_i, p, \rho_i, \nu_i)\stackrel{d_{GH}}{\longrightarrow} (M_{\infty}, p_{\infty},\rho_{\infty}, \nu_{\infty})$ in the measured Gromov-Hausdorff sense, where $\nu_{\infty}$ is the \textbf{renormalized limit measure} defined as in Section $1$ of \cite{CC1}. 


Let $Z$ be a metric space and let $\nu$ be a Borel measure on $Z$. As in Section $2$ of \cite{CC2}, we define the associated Hausdorff measure in codimension $1$ (denoted as $\nu_{-1}$) as follows. Fix $\delta> 0$ and $U\subset Z$, let $\mathcal{B}= \{B_{r_i}(q_i)\}$ be a covering of $U$ with $r_i< \delta$, for all $i$. Put
\begin{align}
(\nu_{-1})_{\delta}(U)= \inf_{\mathcal{B}}\sum_i r_i^{-1} \nu\big(B_{r_i}(q_i)\big) \label{def of HC 0}
\end{align}
and 
\begin{align}
\nu_{-1}(U)= \lim_{\delta\rightarrow 0} (\nu_{-1})_{\delta}(U) \label{def of HC}
\end{align}

\begin{definition}\label{def conic measure}
{On a metric cone $(C(X), dr^2+ r^2dX)$, $\nu$ is called \textbf{conic measure of power $\kappa$},  and $\kappa$ is a positive constant denoted as $\mathtt{p}(\nu)$, if for any $\Omega\subset\subset C(X)$, 
\begin{align}
\nu (\Omega)= \int_0^{\infty} r^{\kappa- 1}dr\int_{X} \chi(\Omega_r) d\nu_{-1} \label{equ sinm}
\end{align}
where $\Omega_r= \{z| z\in \Omega, r(z)= r\}$, $\chi(\cdot)$ is the characteristic function on $C(X)$.
}
\end{definition}

If $(M^n, g)$ is a complete manifold with $Rc\geq 0$ and maximal volume growth, from Theorem \ref{thm metric cone at infinity} above and Theorem $5.9$ in \cite{CC1}, every tangent cone at infinity of $M^n$ is a metric cone, with the unique corresponding renormalized limit measure, which is a conic measure of power $n$. In collapsing case, our definition of conic measure will play the role of co-area formula on metric cones in non-collapsing case, which was showed in Section $7$ of \cite{Honda}.

Assume that $(M^n, g)$ is a complete manifold with $Rc\geq 0$, all tangent cones at infinity are metric cones and every renormalized limit measure is conic measure, we define the set of all tangent cones at infinity of $M^n$ with renormalized limit measure as $\mathscr{M}(M)\vcentcolon= \{(C(X), \nu)|\ C(X)$ is the metric tangent cone at infinity of $M^n$, $\nu$ is the conic renormalized limit measure$\}$.

From \cite{CC3} (also see \cite{Cheeger}), there exists a self-adjoint Laplace operator $\Delta_{(C(X), \nu)}$ on $(C(X), \nu)\in \mathscr{M}(M)$. From (\ref{def of HC 0}) and (\ref{def of HC}), $\nu$ induces a natural measure $\nu_{-1}$ on $X$, which satisfies a volume doubling property. Similar argument as in \cite{Ding} (see Section $4$ there),  weak Poincar\'e inequality also holds on $(X, \nu_{-1})$. Hence from \cite{CC3} (also see \cite{Cheeger}),  volume doubling property, weak Poincar\'e inequality and the rectifiability of the cross section $X$ yields the existence of a self-adjoint positive Laplace operator $\Delta_{(X, \nu_{-1})}$ on $(X, \nu_{-1})$. 

When $\mathcal{H}^1(X)> 0$ where $\mathcal{H}^i$ is $i$-dimensional Hausdorff measure, $L^2(X)$ is an infinite dimensional Hilbert space. Now from Rellich-type Compactness Theorem (Theorem $4.9$ of \cite{Honda1}, also see the Appendix of \cite{Xu}), similar as the standard elliptic theory on compact manifolds (see Chapter $6$ in \cite{Warner} etc.), on compact metric measure space $(X, \nu_{-1})$, we have an orthonormal basis $\{\varphi_i(x)\}_{i= 1}^{\infty}$ for $L^2(X)$, and a sequence $0=\lambda_1< \lambda_2\leq \lambda_3\leq \cdots$, $\lim_{i\rightarrow \infty} \lambda_i= \infty$, such that 
\begin{align}
\Delta_{(X, \nu_{-1})} \varphi_i(x)= -\lambda_i \varphi_i(x) \label{eigen func of lap}
\end{align}

Now assume that $(M^n, g)$ is a complete manifold with $Rc\geq 0$, all tangent cones at infinity are metric cones and every renormalized limit measure is conic measure, and $\mathcal{H}^{1}(X)> 0$. Then we have the following proposition:
\begin{prop}\label{prop A. 2 Laplacian's formula}
{For any $(C(X), \nu)\in \mathscr{M}(M)$, assume that conic measure $\nu$ is of power $\kappa> 0$, for any $u(x)\in H_0^1(C(X))$, 
\begin{align}
\Delta_{(C(X), \nu)} u= \frac{\partial^2 u}{\partial r^2}+ \frac{\kappa- 1}{r}\frac{\partial u}{\partial r}+ \frac{1}{r^2}\Delta_{(X, \nu_{-1})} u \label{A. 2 Laplacian's formula}
\end{align}
}
\end{prop}

\pf
{For any $w\in H_0^1 (C(X))$, by definition of $\Delta_{(C(X), \nu)}$ and $\Delta_{(X, \nu_{-1})}$ (see Section $6$ of \cite{CC3}), we can use integration by parts, combining with the definition of conic measure (\ref{equ sinm}), then
\begin{align}
&\int_{C(X)} \Delta_{(C(X), \nu)} u\cdot w d\nu= -\int_{C(X)} \nabla u\cdot \nabla w d\nu= -\int_{0}^{\infty} r^{\kappa- 1} \Big(\int_X \nabla u\cdot \nabla w d\nu_{-1}\Big) dr \nonumber \\
&\quad \quad =  -\int_{0}^{\infty} r^{\kappa- 1} \Big(\int_X (\nabla_ru+ \frac{1}{r}\nabla_x u)\cdot (\nabla_r w+ \frac{1}{r}\nabla_x w) d\nu_{-1}(x)\Big) dr \label{sepa vari}\\
&\quad\quad = -\int_0^{\infty} r^{\kappa- 3} \Big(\int_X \nabla_x u\cdot \nabla_x w d\nu_{-1}(x)\Big)dr- \int_0^{\infty} r^{\kappa- 1} \int_X \nabla_r u\cdot \nabla_r w d\nu_{-1}(x) \nonumber \\
&\quad \quad = \int_0^{\infty} r^{\kappa- 1} \int_X \frac{1}{r^2}\Delta_{(X, \nu_{-1})} u \cdot w + \int_{0}^{\infty} r^{\kappa- 1} \int_X \Big(\frac{\partial^2 u}{\partial r^2}+ \frac{\kappa- 1}{r}\frac{\partial u}{\partial r}\Big) \cdot w \nonumber \\
&\quad \quad = \int_{C(X)} \Big[\frac{\partial^2 u}{\partial r^2}+ \frac{\kappa- 1}{r}\frac{\partial u}{\partial r}+ \frac{1}{r^2}\Delta_{(X, \nu_{-1})} u\Big] \cdot w d\nu \label{lap 1}
\end{align}
where (\ref{sepa vari}) follows from the metric cone structure of $C(X)$.

From (\ref{lap 1}), we obtain (\ref{A. 2 Laplacian's formula}).
}
\qed

The following corollary is similar as Theorem $1.11$ of \cite{CM-JDG-97} (also see \cite{Cheeger79}), for completeness we provide its proof here following the argument in \cite{CM-JDG-97}.
\begin{cor}\label{cor expression of harmonic function}
{If $u$ is a harmonic function on $(C(X), \nu)$ with respect to $\Delta_{(C(X), \nu)}$, then 
\begin{align}
u(r, x)= \sum_{i= 1}^{\infty} c_i r^{\alpha_i}\varphi_i(x) \label{expression of harmonic function}
\end{align}
where $c_i$, $\alpha_j\geq 0$ are constants, and $\varphi_j$, $\lambda_j= \alpha_j\big(\kappa+ \alpha_j- 2\big)$ are defined in (\ref{eigen func of lap}).
}
\end{cor}

\pf
{We can assume that $u(0)= 0$. By the spectral theorem applied on $(X, \nu_{-1})$, 
\begin{align}
u(1, x)= \sum_{j= 0}^{\infty} a_j \varphi_j(x) \label{ehf 1}
\end{align}
where the convergence is in $L^2(X, \nu_{-1})$ sense.

On the other hand, from Proposition \ref{prop A. 2 Laplacian's formula}, it is not hard to prove that 
\begin{align}
\hat{u}(r, x)= \sum_{j= 0}^{\infty} a_j r^{\alpha_j}\varphi_j(x) \nonumber
\end{align}
is a harmonic function on $(C(X), \nu)$, where $\alpha_j(\alpha_j+ \kappa- 2)= \lambda_j$ and $\alpha_j\geq 0$.

Now consider the harmonic function 
\begin{align}
\tilde{u}(r, x)= u(r, x)- \hat{u}(r, x) \nonumber
\end{align}
From (\ref{ehf 1}), $\tilde{u}$ vanishes on $\partial B_1\subset C(X)$ and at the vertex $0$. Then by the maximum principle, $\tilde{u}\equiv 0$. Hence (\ref{expression of harmonic function}) follows.

}
\qed

And we also define $\mathscr{S}(M)$ the spectrum at infinity of $(M^n, g)$ and $\mathscr{D}(M)$ the degree spectrum at infinity of $(M^n, g)$:
\begin{align}
\mathscr{S}(M)&\vcentcolon= \{\lambda|\ \lambda= \lambda_j(X, \nu_{-1})\ for \ some\ positive\ interger\  j\ and \ (C(X), \nu)\in \mathscr{M}(M)\} \nonumber \\
\mathscr{D}(M)&\vcentcolon= \{\alpha\geq 0|\ \alpha\big(\kappa+ \alpha- 2\big)= \lambda\ for \ some\ \lambda= \lambda_j(X, \nu_{-1})\in \mathscr{S}(M)\  and \ \kappa= \mathtt{p}(\nu)\} \nonumber
\end{align}

We also define the convergence concept for functions on manifolds $\{M_i^n\}$ as the following, it is called ``uniform convergence in Gromov-Hausdorff topology", for simplification, sometimes it is written as "uniform convergence in G-H topology".

\begin{definition}[Uniform Convergence in G-H topology]\label{def uniform convergence}
{Suppose 
\begin{align}
K_i\subset M_i^n\stackrel{d_{GH}}{\longrightarrow} K_{\infty}\subset M_{\infty} \nonumber
\end{align}
Assume that $\{f_i\}_{i= 1}^{\infty}$ are functions on $M_i^n$, $f_{\infty}$ is a function on $M_{\infty}$. and $\Phi_{i}: K_{\infty}\rightarrow K_i$ are $\epsilon_i$-Gromov-Hausdorff approximations, $\lim_{i\rightarrow \infty}\epsilon_i= 0$. If $f_i\circ \Phi_i$ converge to $f_{\infty}$ uniformly, we say that $f_i\rightarrow f_{\infty}$ uniformly over $K_i\stackrel{d_{GH}}{\longrightarrow} K_{\infty}$.
}
\end{definition}

In the rest of this section, unless explicitly stated, $(M^n, g)$ is an $n$-dimensional complete manifold with $Rc\geq 0$ and maximal volume growth. We restrict our discussion to the case of $n\geq 3$, fix $p\in M^n$, let $G(x)$ denote the minimal positive Green's function on $M^n$ with singularity at $p$. And as in \cite{CM}, we will normalize $G(x)$ by 
\begin{align}
\Delta G(x)= (2- n)V_1^{n -1}(\pi)\delta_p(x) \label{normalization of G(x)}
\end{align}
From \cite{Var} (also see \cite{Li}) and the maximal volume growth of the manifold, we know that $G(x)$ exists. Set 
\begin{align}
b(x)= \Big(\frac{V_M}{V_0^n(1)}G(x)\Big)^{\frac{1}{2- n}}\ , \quad\quad\quad \rho(x)= d(p, x)\label{def of b(x)}
\end{align}
Note when $M^n$ is $\mathbb{R}^n$, the function $b(x)$ is just the distance function $\rho(x)$. We also use $B(r)$ to denote the geodesic ball centered at $p$ with radius $r$ on $M$. And we have the following fact:
\begin{align}
\lim_{\rho(x)\rightarrow 0} \frac{b(x)}{\rho(x)}= \Big(\frac{V_M}{V_0^n(1)}\Big)^{\frac{1}{2- n}} \label{estimate of b(0)}
\end{align}

We collect some important facts about $b(x)$ proved by Cheeger and Colding \cite{CC}, Colding and Minicozzi \cite{CM}, Colding \cite{Colding} in the following.
\begin{theorem}[\cite{CC}, \cite{CM}, \cite{Colding}]\label{thm collecting facts}
\begin{align}
\lim_{r\rightarrow \infty} \frac{\int_{b(x)\leq r} \Big||\nabla b|^2- 1\Big|^2 dx}{\mathrm{Vol} (b(x)\leq r)}&= \lim_{r\rightarrow \infty} \frac{\int_{b(x)\leq r} \Big|Hess(b^2)- 2g\Big|^2 dx}{\mathrm{Vol} (b(x)\leq r)}= 0 \ ;\label{estimate of integral of b} \\
\lim_{\rho(x)\rightarrow \infty} \frac{b(x)}{\rho(x)}&= 1\ ;  \quad \quad \quad |\nabla b|\leq 1 \label{estimate of b}
\end{align}
where $g$ is the metric tensor on $M^n$.
\end{theorem}

Let us recall the definition of frequency function in \cite{CM-JDG-97}, we firstly define:
\begin{align}
I_{u}(r)&= r^{1- n}\int_{b(x)= r} u^2|\nabla b| dx       \label{def of I(r)} \\
D_{u}(r)&= r^{2- n}\int_{b(x)\leq r} |\nabla u|^2 dx \ ,  \quad \quad F_{u}(r)= r^{3- n} \int_{b(x)= r} \Big|\frac{\partial u}{\partial n}\Big|^2 |\nabla b| dx \label{def of D(r)} 
\end{align}
then the frequency function is defined by 
\begin{align}
\mathscr{F}_{u}(r)= \frac{D_{u}(r)}{I_{u}(r)} \label{def of frequency function}
\end{align} 
where $u(x)$ is a harmonic function defined on $\{b(x)\leq r\}$. 

Using the fact that $u$ is harmonic, differentiating (\ref{def of I(r)}), we get 
\begin{align}
I'_{u}(r)= 2\frac{D_{u}(r)}{r}\geq 0 \label{I'(r)}
\end{align}

From (\ref{I'(r)}), $I_1(r)$ is constant. Then by the fact (\ref{estimate of b(0)}), it is not hard to see that 
\begin{align}
I_1(r)= nV_M \label{I_1(r)}
\end{align}

We further define two quantities which are technically easier to be dealt with, comparing with $D_{u}$ and $\mathscr{F}_{u}$.
\begin{align}
E_{u}(r)= r^{2- n}\int_{b(x)\leq r} |\nabla u|^2 |\nabla b|^2 dx \ , \quad \quad 
\mathscr{W}_{u}(r)= \frac{E_{u}(r)}{I_{u}(r)} \label{def of W(r)}
\end{align}

Sometimes for simplification, we omit the subscript $u$ in $I_{u}(r)$, $\cdots$, $\mathscr{W}_{u}(r)$ when the context is clear, and use $I(r)$, $\cdots$, $\mathscr{W}(r)$ instead.

When $M^n$ is a complete manifold with $Rc\geq 0$ and maximal volume growth, $r_j\rightarrow \infty$, assume that the rescaled manifolds $(M^n, p, r_j^{-2}g)$ converge to some tangent cone at infinity $M_{\infty}$ in the Gromov-Hausdorff sense,  From Theorem $0.1$ of \cite{CM}, and Theorem $3.21$, Corollary $4.22$ of \cite{Ding}, we have the following proposition:

\begin{prop}\label{prop b is continuous under G-H convergence}
{If $K_j$ and $K_{\infty}$ are compact subsets of $(M^n, p, r_j^{-2}g)$ and $M_{\infty}$ respectively, suppose $K_j\stackrel{d_{GH}}{\longrightarrow} K_{\infty}$, then $b_j\rightarrow b_{\infty}$ uniformly over $K_j\stackrel{d_{GH}}{\longrightarrow} K_{\infty}$, where $b_j$ and $b_{\infty}$ are defined by (\ref{def of b(x)}) on $(M^n, p, r_j^{-2}g)$ and $(M_{\infty}, p_{\infty},\rho_{\infty})$ respectively, furthermore $b_{\infty}= \rho_{\infty}$.
}
\end{prop}

\section{Three Circles Theorems for harmonic functions}\label{SECTION 3}
Different types of Three Circles Theorems were proved by Simon \cite{Simon}, Cheeger and Tian \cite{CT}, Colding, DeLellis and Minicozzi \cite{CDM} in different contexts. Also see Zhang \cite{Zhang}, Ding \cite{Ding2} for harmonic functions on manifolds and Liu \cite{Liu} for holomorphic functions on K\"ahler manifolds. 

However, to study Question \ref{ques 1.3}, Conjecture \ref{conj 5.4} and Conjecture \ref{conj 1.2}, we need to do some modification to get the Three Circles Theorems applicable on those problems.

\begin{lemma}\label{lemma 3.1 induction}
{For $\{w_i\}_{i= 1}^{\infty}$, $w_i\geq 0$, if 
\begin{align}
\sum_{i= 1}^{\infty} w_i\leq \sum_{i= 1}^{\infty} 2^{2(\alpha- \alpha_i)} w_i \label{3.1.1}
\end{align}
then
\begin{align}
\sum_{i= 1}^{\infty} 2^{-2\alpha_i}w_i \leq  \sum_{i= 1}^{\infty} 2^{2(\alpha- 2\alpha_i)} w_i \label{3.1.2} 
\end{align}
where $0= \alpha_1< \alpha_2\leq \alpha_3\leq \cdots$, and $\alpha> 0$. Furthermore, the equality in (\ref{3.1.2}) holds if and only if $w_i= 0$ for all $i$ satisfying $\alpha_i\neq \alpha$. 
}
\end{lemma}

\pf
{(\ref{3.1.1}) is equivalent to 
\begin{align}
\sum_{\alpha_i\neq \alpha_1} w_i(1- 2^{2(\alpha- \alpha_i)})\leq w_1(2^{2\alpha}- 1)  \label{3.1.3}
\end{align}
and (\ref{3.1.2}) is equivalent to 
\begin{align}
\sum_{\alpha_i\neq \alpha_1} 2^{-2\alpha_i}w_i(1- 2^{2(\alpha- \alpha_i)})\leq  w_1(2^{2\alpha}- 1)  \label{3.1.4}
\end{align}

Note 
\begin{align}
\sum_{\alpha_i\neq \alpha_1} 2^{-2\alpha_i}w_i(1- 2^{2(\alpha- \alpha_i)})\leq \sum_{\alpha_i\neq \alpha_1} 2^{-2\alpha}w_i(1- 2^{2(\alpha- \alpha_i)}) \label{3.1.7}
\end{align}

If $\sum_{\alpha_i\neq \alpha_1}w_i(1- 2^{2(\alpha- \alpha_i)})\leq 0$, from (\ref{3.1.7}), 
\begin{align}
\sum_{\alpha_i\neq \alpha_1} 2^{-2\alpha_i} w_i(1- 2^{2(\alpha- \alpha_i)})\leq 0\leq w_1(2^{2\alpha}- 1)  \nonumber
\end{align}
If $\sum_{\alpha_i\neq \alpha_1}w_i(1- 2^{2(\alpha- \alpha_i)})> 0$, from (\ref{3.1.3})
\begin{align}
\sum_{\alpha_i\neq \alpha_1} 2^{-2\alpha_i} w_i(1- 2^{2(\alpha- \alpha_i)})\leq 2^{-2\alpha}w_1(2^{2\alpha}- 1)< w_1(2^{2\alpha}- 1)   \nonumber 
\end{align} 

Hence (\ref{3.1.4}) is proved, and (\ref{3.1.2}) is obtained. Check the above argument carefully, it is easy to find that the equality in (\ref{3.1.2}) holds if and only if $w_i= 0$ for all $i$ satisfying $\alpha_i\neq \alpha$.
}
\qed

On $(M^n, g, \mu)$ where $\mu$ is a Borel regular measure on $M^n$, define the $J$-function of $u$ as the following:
\begin{align}
J_u(r)= \frac{1}{\mu(B(r))}\int_{B(r)} u^2 d\mu \label{def of J}
\end{align}
Unless otherwise mentioned, for $J_u(r)$ in (\ref{def of J}), the measure $\mu$ will be assumed as the volume measure determined by the metric $g$.

\begin{theorem}\label{thm 3.6 induction of J}
{Let $(M^n, g)$ be a complete manifold with nonnegative Ricci curvature, assume that every tangent cone at infinity of $M^n$ is a metric cone $C(X)$ with the conic renormalized limit measure of power $\kappa\geq 2$, and $\mathcal{H}^{1}(X)> 0$. If $\alpha\notin \mathscr{D}(M)$, then there exists integer $k_0= k_0(\alpha)> 1$, such that for $r\geq k_0$, and $u(x)$ harmonic over $B(r)\subset (M^n, g)$, 
\begin{equation}\label{3.6.1}
J_u(r)\leq 2^{2\alpha} J_u(\frac{r}{2})
\end{equation}
implies 
\begin{align}
J_u(\frac{r}{2})\leq 2^{2\alpha} J_u(\frac{r}{4})   \label{3.6.2}
\end{align}
}
\end{theorem}

\pf
{By contradiction. If Theorem \ref{thm 3.6 induction of J} is not true, then there exists a sequence $\{r_l\}$, $r_l\rightarrow \infty$, and the corresponding harmonic functions $u_{l}$ such that the following inequalities hold:
\begin{align}
J_{u_l}(r_l)\leq 2^{2\alpha} J_{u_l}(\frac{r_l}{2})\ , \quad \quad J_{u_l}(\frac{r_l}{2})> 2^{2\alpha} J_{u_l}(\frac{r_l}{4}) \label{3.6.3}
\end{align}
Using the assumptions about tangent cones with renormalized limit measure, combining the knowledge about measured Gromov-Hausdorff convergence, without loss of generality (by choosing subsequence of $\{u_l\}$), we assume $(M^n, p, \rho_l, \nu_i)$ converges to $(C(Y), p_{\infty}, \rho_{\infty}, \nu_{\infty})$ in the measured Gromov-Hausdorff sense as in $\textbf{Blow Down Setup}$ of Section \ref{SECTION 2}, and $(C(Y), \nu_{\infty})$ is a metric cone with conic measure.

Clearly (\ref{3.6.3}) implies $J_{u_l}(\frac{r_l}{2})\neq 0$, define 
\begin{align}
\tilde{u}_l= \frac{u_l}{\Big(J_{u_l}(\frac{r_l}{2})\Big)^{\frac{1}{2}}} \nonumber
\end{align}

Look at $\tilde{u}_l$ as the function on $B_l(1)\subset (M^n, g_l)$, from (\ref{3.6.3}) 
\begin{align}
J_{\tilde{u}_l}^{(l)}(1)\leq 2^{2\alpha} J_{\tilde{u}_l}^{(l)}(\frac{1}{2}) \ , \quad \quad J_{\tilde{u}_l}^{(l)}(\frac{1}{2})> 2^{2\alpha} J_{\tilde{u}_l}^{(l)}(\frac{1}{4})\label{3.6.6} 
\end{align}
where $J_{\tilde{u}_l}^{(l)}$ is the $J$-function of $\tilde{u}_l$ on manifold $(M^n, g_l, \nu_l)$. Also we have
\begin{align}
J_{\tilde{u}_l}^{(l)}(\frac{1}{2})= J_{\tilde{u}_l}(\frac{r_l}{2})= 1 \label{3.6.7}
\end{align}

From Theorem $1.2$ in \cite{LS} and Cheng-Yau's gradient estimate in \cite{CY}, we have the following estimates:
\begin{align}
&\sup_{B^{(l)}(1-\theta)} |\tilde{u}_l|\leq C(n, p, \theta)\Big[J_{\tilde{u}_l}^{(l)}(1)\Big]^{\frac{1}{2}}\leq C(n, p, \theta, \alpha)\Big[J_{\tilde{u}_l}^{(l)}(\frac{1}{2})\Big]^{\frac{1}{2}}= C(n, p, \theta, \alpha)\nonumber \\
&\sup_{B^{(l)}(1-\theta)} |\nabla \tilde{u}_l|_{g_l}\leq  C(n, p, \theta, \alpha) \nonumber
\end{align}

So for any $\theta\in (0, 1)$, $\tilde{u}_l$ and $|\nabla \tilde{u}_l|$ are uniformly bounded over $B^{(l)}(1-\theta)$. By Harnack's convergence theorem in the Gromov-Hausdorff sense (see \cite{Ding}, also \cite{Xu}), we get that $\{\tilde{u}_l\}$ converges uniformly on compact subsets of $B_{\infty}(1)\subset C(Y)$ to $u_{\infty}$, and $u_{\infty}$ is harmonic over $B_{\infty}(1)$. Hence,
\begin{align}
\lim_{l\rightarrow \infty} J_{\tilde{u}_l}^{(l)}(\frac{1}{2})= J_{u_{\infty}}(\frac{1}{2}) \label{3.6.8}
\end{align}
where $J_{u_{\infty}}$ is the $J$-function of $u_{\infty}$ on $(C(Y), \rho_{\infty}, \nu_{\infty})$ defined as in (\ref{def of J}).

By $u_{\infty}$ is harmonic over $B_{\infty}(1)\subset C(Y)$, as in (\ref{expression of harmonic function}) we can write 
\begin{align}
u_{\infty}= \sum_{i= 1}^{\infty} c_i r^{\alpha_i} \varphi_i(x) \label{3.6.9}
\end{align}
where $\{\varphi_i(x)\}$ are the eigenfunctions of $\Delta_Y$ on $Y$, also the orthonormal basis for $L^2(Y)$, $\Delta_Y \varphi_i(x)= -\lambda_i \varphi_i(x)$, $\lambda_i= \alpha_i(\alpha_i+ n- 2)$ and $\alpha_i\geq 0$.

From (\ref{3.6.6}), we get
\begin{align}
\int_{B_{\infty}(1)} u_{\infty}^2 d\nu_{\infty}&= \lim_{k\rightarrow \infty} \int_{B_{\infty}\big(1- \frac{1}{k}\big)} u_{\infty}^2 d\nu_{\infty}= \lim_{k\rightarrow \infty} \lim_{l\rightarrow \infty} \int_{B_{l}\big(1- \frac{1}{k}\big)} \tilde{u}_l^2 d\nu_{l}\nonumber \\
&\leq \lim_{l\rightarrow \infty} \int_{B_{l}(1)} \tilde{u}_l^2 d\nu_{l}\leq 2^{2\alpha} \lim_{l\rightarrow \infty} \frac{V(r_l)}{V(\frac{r_l}{2})}\int_{B_{l}(\frac{1}{2})} \tilde{u}_l^2 d\nu_l \nonumber \\
&= 2^{2\alpha+ \kappa} \int_{B_{\infty}(\frac{1}{2})} u_{\infty}^2 d\nu_{\infty}
\label{3.6.10}\\
\int_{B_{\infty}(\frac{1}{2})} u_{\infty}^2 d\nu_{\infty} &\geq 2^{2\alpha+ \kappa}  \int_{B_{\infty}(\frac{1}{4})} u_{\infty}^2 d\nu_{\infty}   \label{3.6.11}
\end{align}
in the last equality of (\ref{3.6.10}) we used the assumption that the renormalized limit measure is conic measure of degree $\kappa$.

Plug (\ref{3.6.9}) into (\ref{3.6.10}) and (\ref{3.6.11}), we get 
\begin{align}
\sum_{i=1}^{\infty} w_i\leq \sum_{i=1}^{\infty} 2^{2(\alpha- \alpha_i)} w_i \ , \quad\quad 
\sum_{i=1}^{\infty} 2^{-2\alpha_i}w_i \geq \sum_{i=1}^{\infty} 2^{2(\alpha- 2\alpha_i)} w_i \nonumber
\end{align}
where $w_i= \frac{c_i^2}{2\alpha_i+ \kappa}$.

From the above two inequalities, by Lemma \ref{lemma 3.1 induction} and the assumption $\alpha\notin \mathscr{D}(M)$, we get $w_i= 0$ and $c_i= 0$, hence $u_{\infty}\equiv 0$. Taking limit in (\ref{3.6.7}), combining (\ref{3.6.8}) and $u_{\infty}= 0$, we obtain
\begin{align}
1=\lim_{l\rightarrow \infty} J_{\tilde{u}_l}^{(l)}(\frac{1}{2})= J_{u_{\infty}}(\frac{1}{2})= 0 \nonumber
\end{align}

It is the contradiction, hence the conclusion is proved.
}
\qed

Recall we defined $I_u$ in (\ref{def of I(r)}) for harmonic functions $u(x)$, we have the other Three Circles Theorems for $I_u$, which will be useful for estimating the frequency of $u(x)$. Before proving the theorem, we firstly need to control the $C^0$ and $C^1$ norm of $u(x)$ by $I_u$, which is achieved by the following lemma.

\begin{lemma}\label{lem 3.4 estimate of u}
{Assume that $(M^n, g)$ is an $n$-dimensional complete manifold with $Rc\geq 0$ and maximal volume growth, $u(p)= 0$ and $u(x)$ is harmonic on $\{b\leq r\}\subset M^n$. Then for any $\theta\in (0, 1)$ and $r> 0$, 
\begin{align}
&\sup_{b\leq (1- \theta)r} |u|^2 \leq C(n, p, \theta, V_M)I_u(r) \label{estimate 1} \\
&\sup_{b\leq (1- \theta)r} |\nabla u|^2 \leq C(n, p, \theta, V_M)r^{-2}I_u(r)  \label{estimate 2}
\end{align}
}
\end{lemma}

\pf
{From $I'(r)= \frac{2D(r)}{r}$, we have
\begin{align}
2\int_{(1- \frac{\theta}{2})r}^{r} \frac{D(s)}{s}= I(r)- I\Big(\big(1- \frac{\theta}{2}\big)r\Big)\leq I(r) \nonumber
\end{align}
hence 
\begin{align}
2\int_{(1- \frac{\theta}{2})r}^{r} s^{n- 2}D(s) ds\leq r ^{n- 1}I(r) \label{3.2.1}
\end{align}
by the definition of $D(r)$, $s^{n- 2}D(s)$ is nondecreasing, and therefore (\ref{3.2.1}) yields 
\begin{align}
\theta r\Big(\big(1- \frac{\theta}{2})r\Big)^{n- 2} D\Big(\big(1- \frac{\theta}{2})r\Big)\leq r^{n- 1}I(r) \nonumber
\end{align}
After simplification, we get 
\begin{align}
D\Big(\big(1- \frac{\theta}{2})r\Big)\leq C(n, \theta)I(r) \label{3.2.2}
\end{align}

Assume that $|\nabla u|^2(x_0)= \sup_{b\leq (1- \theta)r} |\nabla u|^2$ and $b(x_0)= r_0\leq (1- \theta)r$. From (\ref{estimate of b}), there exists $r_1= C(p, V_M, \theta)r> 0$ such that $B_{x_0}(r_1)\subset \{b\leq \big(1- \frac{\theta}{2}\big)r\}$, hence by Theorem $1.2$ in \cite{LS}, 
\begin{align}
|\nabla u|^2(x_0)\leq \frac{C(n)}{V\Big(B_{x_0}(r_1)\Big)} \int_{B_{x_0}(r_1)} |\nabla u|^2 \leq C(n, p, \theta, V_M)r^{-2} D\Big(\big(1- \frac{\theta}{2}\big)r\Big) \nonumber
\end{align}

Combining (\ref{3.2.2}), we obtain
\begin{align}
\sup_{b\leq (1- \theta)r} |\nabla u|^2\leq  C(n, p, \theta, V_M)r^{-2} I(r) \label{3.2.3}
\end{align} 

By integrating (\ref{3.2.3}) along geodesics starting at $p$ and using $u(p)= 0$, 
\begin{align}
\sup_{b\leq (1- \theta)r} |u|^2\leq C(n, p, \theta, V_M) I(r) \nonumber
\end{align}
}
\qed

\begin{theorem}\label{thm 3.5 induction of I}
{Assume that $(M^n, g)$ is an $n$-dimensional complete manifold with $Rc\geq 0$ and maximal volume growth, $\alpha\notin \mathscr{D}(M)$, then there exists integer $k_0= k_0(\alpha)> 1$, such that for $r\geq k_0$, $u(x)$ harmonic over $B(r)\subset (M^n, g)$ and $u(p)= 0$, 
\begin{equation}\label{3.5.1}
I_u(r)\leq 2^{2\alpha} I_u(\frac{r}{2})
\end{equation}
implies 
\begin{align}
I_u(\frac{r}{2})\leq 2^{2\alpha} I_u(\frac{r}{4})   \label{3.5.2}
\end{align}
}
\end{theorem}

\pf
{By contradiction. If Theorem \ref{thm 3.5 induction of I} is not true, then there exists a sequence $\{r_l\}$, $r_l\rightarrow \infty$, and the corresponding harmonic functions $u_{l}$ such that the following inequalities hold:
\begin{align}
I_{u_l}(r_l)\leq 2^{2\alpha} I_{u_l}(\frac{r_l}{2}) \ , \quad \quad I_{u_l}(\frac{r_l}{2})> 2^{2\alpha} I_{u_l}(\frac{r_l}{4})\ , \quad \quad u_l(p)= 0 \label{3.5.4}
\end{align}
Without loss of generality (by choosing subsequence of $\{u_l\}$), we can assume that
\begin{align}
(M^n, p, \rho_l)\stackrel{d_{GH}}{\longrightarrow} (C(X), p_{\infty}, \rho_{\infty}) \label{convergence}
\end{align}
where $\rho_l= g_l= r_l^{-2}g$ is the rescaled metric, and $C(X)$ is one tangent cone at infinity of $(M^n, g)$, which is a metric cone by Theorem \ref{thm metric cone at infinity}.

Clearly (\ref{3.5.4}) implies $I_{u_l}(\frac{r_l}{2})\neq 0$, define 
\begin{align}
\tilde{u}_l= \frac{u_l}{\Big(I_{u_l}(\frac{r_l}{2})\Big)^{\frac{1}{2}}} \nonumber
\end{align}

Look at $\tilde{u}_l$ as the function on $B_l(1)\subset (M^n, g_l)$, we have
\begin{align}
I_{\tilde{u}_l}^{(l)}(1)\leq 2^{2\alpha} I_{\tilde{u}_l}^{(l)}(\frac{1}{2}) \ , \quad \quad 
I_{\tilde{u}_l}^{(l)}(\frac{1}{2})> 2^{2\alpha} I_{\tilde{u}_l}^{(l)}(\frac{1}{4})\ , \quad \quad \tilde{u}_l(p)= 0\label{3.5.6} 
\end{align}
where $I_{\tilde{u}_l}^{(l)}$ is the frequency function of $\tilde{u}_l$ on manifold $(M^n, g_l)$, and 
\begin{align}
I_{\tilde{u}_l}^{(l)}(\frac{1}{2})= I_{\tilde{u}_l}(\frac{r_l}{2})= 1 \label{3.5.7}
\end{align}

From Lemma \ref{lem 3.4 estimate of u} and (\ref{3.5.6}), we have the following estimates,
\begin{align}
\sup_{b_l\leq (1-\theta)} |\tilde{u}_l|&\leq C(n, p, \theta, V_M)\Big[I_{\tilde{u}_l}^{(l)}(1)\Big]^{\frac{1}{2}}\leq C(n, p, \theta, V_M, \alpha)\Big[I_{\tilde{u}_l}^{(l)}(\frac{1}{2})\Big]^{\frac{1}{2}} \nonumber \\
&= C(n, p, \theta, V_M, \alpha)\nonumber \\
\sup_{b_l\leq (1-\theta)} |\nabla \tilde{u}_l|_{g_l}&\leq  C(n, p, \theta, V_M, \alpha) \nonumber
\end{align}
where $b_l$ is the $b(x)$ function defined as in (\ref{def of b(x)}) on $(M^n, p, g_l)$.

So for any $\theta\in (0, 1)$, $\tilde{u}_l$ and $|\nabla \tilde{u}_l|$ are uniformly bounded over $\{b_l\leq (1-\theta)\}$. By Harnack's convergence theorem in the Gromov-Hausdorff sense (see \cite{Ding}, also \cite{Xu}), we get that $\{\tilde{u}_l\}$ converges uniformly on compact subsets of $\{b_{\infty}< 1\}\subset C(X)$ to $u_{\infty}$, and $u_{\infty}$ is harmonic over $\{b_{\infty}< 1\}$. From Proposition $3.4$ in \cite{Honda}, we have
\begin{align}
\lim_{l\rightarrow \infty} I_{\tilde{u}_l}^{(l)}(\frac{1}{2})= I_{u_{\infty}}(\frac{1}{2}) \label{3.5.8}
\end{align}

By $u_{\infty}$ is harmonic over $\{b_{\infty}< 1\}$, as in (\ref{expression of harmonic function}) we can write
\begin{align}
u_{\infty}= \sum_{i= 1}^{\infty} c_i r^{\alpha_i} \varphi_i(x) \label{3.5.9}
\end{align}
where $\{\varphi_i(x)\}$ are the eigenfunctions of $\Delta_X$ on $X$, also the orthonormal basis for $L^2(X)$, $\Delta_X \varphi_i(x)= -\lambda_i \varphi_i(x)$, $\lambda_i= \alpha_i(\alpha_i+ n- 2)$ and $\alpha_i\geq 0$.

From (\ref{3.5.6}), note $b_{\infty}= \rho_{\infty}$ on $C(X)$, again by Proposition $3.4$ in \cite{Honda} 
\begin{align}
\int_{b_{\infty}= 1} u_{\infty}^2 &= \lim_{k\rightarrow \infty} \int_{b_{\infty}= \big(1- \frac{1}{k}\big)} u_{\infty}^2= \lim_{k\rightarrow \infty} \Big\{\Big(1- \frac{1}{k}\Big)^{n- 1}\lim_{l\rightarrow \infty} I_{\tilde{u}_l}^{(l)}\big(1- \frac{1}{k}\big) \Big\}\nonumber \\
&\leq \lim_{l\rightarrow \infty} I_{\tilde{u}_l}^{(l)}(1)\leq \lim_{l\rightarrow \infty} 2^{2\alpha} I_{\tilde{u}_l}^{(l)}(\frac{1}{2})= 2^{2\alpha+ n- 1} \int_{b_{\infty}= \frac{1}{2}} u_{\infty}^2
\label{3.5.10}\\
\int_{b_{\infty}= \frac{1}{2}} u_{\infty}^2 &\geq 2^{2\alpha+ n- 1} \int_{b_{\infty}= \frac{1}{4}} u_{\infty}^2   \label{3.5.11}
\end{align}
in the first inequality of (\ref{3.5.10}) we used the fact $I(r)$ is nondecreasing in $r$.

Plug (\ref{3.5.9}) into (\ref{3.5.10}) and (\ref{3.5.11}), we get 
\begin{align}
\sum_{i=1}^{\infty} w_i\leq \sum_{i=1}^{\infty} 2^{2(\alpha- \alpha_i)} w_i \ , \quad \quad  \sum_{i=1}^{\infty} 2^{-2\alpha_i}w_i \geq \sum_{i=1}^{\infty} 2^{2(\alpha- 2\alpha_i)} w_i \nonumber
\end{align}
where $w_i= c_i^2$. From the above two inequalities and Lemma \ref{lemma 3.1 induction}, $c_i= 0$, hence $u_{\infty}\equiv 0$.

Taking limit in (\ref{3.5.7}), combining (\ref{3.5.8}) and $u_{\infty}= 0$, we obtain
\begin{align}
1=\lim_{l\rightarrow \infty} I_{\tilde{u}_l}^{(l)}(\frac{1}{2})= I_{u_{\infty}}(\frac{1}{2})= 0 \nonumber
\end{align}

It is the contradiction, hence the conclusion is proved.
}
\qed

\section{The existence of harmonic functions with polynomial growth}\label{SECTION 4}

In the following lemma, we assume that $(M^n, g)$ is a complete manifold with $Rc\geq 0$, and every tangent cone at infinity of $M^n$ with renormalized limit measure is a metric cone $C(X)$ with conic measure of power $\kappa\geq 2$, and $\mathcal{H}^{1}(X)> 0$.

\begin{lemma}\label{lemma 5.1 construction of u_i}
{Assume $u_{\infty}$ is harmonic on $C(X)\in \mathscr{M}(M)$ and $u_{\infty}(p_{\infty})= 0$, then there exist $R_i\rightarrow \infty$, $B(R_i)\subset M^n$, such that $\lim_{i\rightarrow \infty} d_{GH} \big(B_i(1), B_{\infty}(1)\big)= 0$, where $B_i(1)\subset (M^n, R_i^{-2}g)$, $B_{\infty}(1)\subset C(X)$, and $u_i$ harmonic on $B_p(R_i)= B_i(1)$ satisfying the following property:
\begin{align}
\lim_{i\rightarrow \infty} |u_i\circ \Psi_{\infty, i}- u_{\infty}|_{L^{\infty}\big(B_{\infty}(1)\big)}=0\ , \quad\quad\quad u_i(p)= 0 \label{5.1.1}
\end{align}
where $\Psi_{\infty, i}: B_{\infty}(1)\rightarrow  B_i(1)$ is an $\epsilon_i$-Gromov-Hausdorff approximation, and $\lim_{i\rightarrow \infty} \epsilon_i= 0$.
}
\end{lemma}

\begin{remark}
{The above lemma was inspired by Theorem $2.1$ of \cite{Ding2}, however there is a small (but new) restriction on $u_i$ ($u_i(p)= 0$) in our statement, which is crucial in our proof of Theorem \ref{thm general existence of poly growth harmonic function}.
}
\end{remark}

\pf
{It follows from Lemma $10.7$ in \cite{Cheeger}, there exists Lipschitz function $\tilde{u}_i$ defined on $(M^n, R_i^{-2}g)$ such that
\begin{align}
\lim_{i\rightarrow \infty} |\tilde{u}_i\circ \Psi_{\infty, i}- u_{\infty}|_{L^{\infty}\big(\overline{B_{\infty}}(1)\big)}= 0\ ; \quad\quad \textbf{Lip}(\tilde{u}_i)\leq C \label{5.1.2}
\end{align}
where $C$ is some positive constant independent of $i$, and 
\begin{align}
\textbf{Lip}(f)\vcentcolon= \sup_{z\in \overline{B_{\infty}}(1)}\liminf_{r\rightarrow 0}\sup_{d(z, y)= r}\frac{|f(y)- f(z)|}{r} \nonumber
\end{align}

Let $\hat{u}_i$ be the solution of the following Dirichlet problem:
\begin{equation}\label{Lap equ}
\left\{
\begin{array}{rl}
\Delta \hat{u}_i&= 0\ , \quad on \ B_i(1)\\
\hat{u}_i&= \tilde{u}_i\ , \quad on\ \partial B_i(1) 
\end{array} \right.
\end{equation}

By (\ref{5.1.2}) and Lemma $10.7$ in \cite{Cheeger}, $\tilde{u}_i|_{\partial B_i(1)}$ is uniformly bounded. From maximum principle on $B_i(1)$, $\hat{u}_i$ are uniformly bounded. 

Let $x_i\in \partial B_i(1)$, $x_i\rightarrow x_{\infty}\in \partial B_{\infty}(1)$. For any $\epsilon> 0$, from (\ref{5.1.2}) there exists $\delta\in (0, 1)$ such that
\begin{align}
|\tilde{u}_i(x)- \tilde{u}_i(x_i)|< \frac{\epsilon}{2} \quad \quad \quad if \ \ |x- x_i|< \delta \nonumber 
\end{align} 

On the cone $C(X)$ there is a unique ray starting from the pole $p_{\infty}$, passing through $x_{\infty}$. Pick a point $q_{\infty}$ on this ray with $d_{\infty}(p_{\infty}, q_{\infty})> d_{\infty}(p_{\infty}, x_{\infty})$, then 
\begin{align}
d_{\infty}(x, q_{\infty})> d(x_{\infty}, q_{\infty})\ , \quad \quad \quad \forall x\in \{z|\ d_{\infty}(z, x_{\infty})< \delta\}\cap \bar{B}_{\infty}(1) \nonumber
\end{align}
Hence we can choose $q_i\rightarrow q_{\infty}$ such that 
\begin{align}
d_{\rho_i} (q_i, x)\geq d_{\rho_i}(q_i, x_i) \ , \quad \quad \quad \forall x\in \{z|\ d_{\rho_i}(z, x_{i})< \delta\}\cap \bar{B}_{i}(1) \nonumber
\end{align}

Consider 
\begin{align}
w_i(x)= d_{\rho_i}(q_i, x_i)^{2- n}- d_{\rho_i}(q_i, x)^{2- n}  \nonumber
\end{align}
By the Laplacian comparison theorem, $\Delta_i w_i\leq 0$, and it is easy to see that
$w_i\geq 0$ on $\{z|\ d_{\rho_i}(z, x_{i})< \delta\}\cap \bar{B}_{i}(1)$, $w_i(x_i)= 0$. Hence it is the barrier function defined as in Section $2.8$ of \cite{GT}.

Now let $\mathcal{M}= \sup_{i}\sup_{x\in \partial B_i(1)} |\tilde{u}_i(x)|< \infty$, using the fact 
\begin{align}
\lim_{i\rightarrow \infty} w_i(y_i)= d_{\infty}(q_{\infty}, x_{\infty})^{2- n}- d_{\infty}(q_{\infty}, y_{\infty})^{2- n}\quad \quad \quad if\ \  y_i\rightarrow y_{\infty} \label{limit fact}
\end{align}
hence there exists constant $k> 0$, which is independent of $i$, such that when $i> > 1$,
\begin{align}
kw_i(x)\geq 2\mathcal{M} \quad \quad \quad if \ \ |x- x_i|\geq \delta \nonumber 
\end{align}

Then it is easy to check
\begin{align}
\Big( \tilde{u}_i(x_i)+ \frac{\epsilon}{2}+ kw_i(x)- \hat{u}_i(x)\big|_{\partial B_i(1)} \Big)&\geq 0\ , \nonumber \\
\Delta\Big( \tilde{u}_i(x_i)+ \frac{\epsilon}{2}+ kw_i(x)- \hat{u}_i(x) \Big) &\leq 0 \nonumber 
\end{align}

From maximum principle, in $B_i(1)$,
\begin{align}
\tilde{u}_i(x_i)+ \frac{\epsilon}{2}+ kw_i(x)\geq \hat{u}_i(x) \nonumber
\end{align}

Similarly in $B_i(1)$, we have
\begin{align}
\tilde{u}_i(x_i)- \frac{\epsilon}{2}- kw_i(x)\leq \hat{u}_i(x) \nonumber
\end{align}

Hence 
\begin{align}
|\hat{u}_i(x)- \tilde{u}_i(x_i)|\leq \frac{\epsilon}{2}+ kw_i(x) \quad \quad x\in B_i(1) \label{uniform cont near boundary}
\end{align}

Note $k$ is independent of $i$, using the fact (\ref{limit fact}) again, we get $\delta_0> 0$, such that for $x_i\in \partial B_i(1)$, $d_{\rho_i}(x, x_i)\leq \delta_0$ implies $|\hat{u}_i(x)- \tilde{u}_i(x_i)|\leq \epsilon$ for any $i> > 1$. In fact, from (\ref{uniform cont near boundary}), we can obtain that $\hat{u}_i$ are uniformly continuous near boundary of $\overline{B_i}(1)$. Combining with the Cheng-Yau's gradient estimate, $\hat{u}_i$ are uniformly continuous on $\overline{B_i}(1)$. 

From Harnack's convergence theorem in the Gromov-Hausdorff sense (see \cite{Ding}, also \cite{Xu}), $\hat{u}_i$ converges to $w_{\infty}$ on $\overline{B_{\infty}}(1)$, i.e.
\begin{align}
\lim_{i\rightarrow \infty} |\hat{u}_i\circ \Psi_{\infty, i}- w_{\infty}|_{L^{\infty}\big(\overline{B_{\infty}}(1)\big)}= 0 \label{5.1.3}
\end{align}

From (\ref{5.1.2}), (\ref{Lap equ}) and (\ref{5.1.3}), we get that $w_{\infty}|_{\partial B_{\infty}(1)}= u_{\infty}|_{\partial B_{\infty}(1)}$. From maximum principle on $C(X)$, $w_{\infty}= u_{\infty}$ on $\overline{B_{\infty}}(1)$. We get 
\begin{align}
\lim_{i\rightarrow \infty} |\hat{u}_i\circ \Psi_{\infty, i}- u_{\infty}|_{L^{\infty}\big(\overline{B_{\infty}}(1)\big)}= 0 \nonumber
\end{align}

Choose $u_i(x)= \hat{u}_i(x)- \hat{u}_i(p)$, note $\hat{u}_i(p)\rightarrow u_{\infty}(p_{\infty})= 0$, then 
\begin{align}
\lim_{i\rightarrow \infty} |u_i\circ \Psi_{\infty, i}- u_{\infty}|_{L^{\infty}\big(\overline{B_{\infty}}(1)\big)}= 0 \nonumber
\end{align}

The conclusion is obtained.
}
\qed

\begin{theorem}\label{thm general existence of poly growth harmonic function}
{Let $(M^n, g)$ be a complete manifold with nonnegative Ricci curvature, assume that every tangent cone at infinity of $M^n$ with renormalized limit measure is a metric cone $C(X)$ with conic measure of power $\kappa\geq 2$, and $\mathcal{H}^{1}(X)> 0$. If there exists $d\notin \mathscr{D}(M)$ and $d> \inf \{\alpha|\ \alpha\in \mathscr{D}(M), \alpha\neq 0\}$, then $dim\big(\mathscr{H}_d(M)\big)\geq 2$.
}
\end{theorem}

\pf
{By assumption, there exists $\alpha_1\in \mathscr{D}(M)$, $\alpha_1\neq 0$ and $\alpha_1< d$. Hence there is $C(X)\in \mathscr{M}(M)$, and $\varphi_1(x)$ is the eigenfunction on $X$ with respect to eigenvalue $\lambda_1= \alpha_1(\alpha_1+ \kappa- 2)$, $\int_{X} |\varphi_1|^2= 1$. Let $u_{\infty}= r^{\alpha_1} \varphi_1(x)$ in Lemma \ref{lemma 5.1 construction of u_i}, then choose $\{u_i\}$ from Lemma \ref{lemma 5.1 construction of u_i}. We have the following lemma:

\begin{lemma}\label{lemma 5.2 estimate of u_i near infinity}
{For any given positive constant $r_0\in (0, 1)$, there exists $i_0= i_0(d- \alpha_1, r_0)> 0$ such that if $i\geq i_0$, for any $r\in [r_0R_i, R_i]$,
\begin{equation}\label{5.2.1}
J_{u_i}(r)\leq 2^{2d} \cdot J_{u_i}(\frac{r}{2})  
\end{equation}
}
\end{lemma}

\pf
{By contradiction. If the lemma is not true, without loss of generality, we can assume that for some $r_0\in (0, 1)$, there exists a subsequence of $\{i\}_{1}^{\infty}$, for simplicity also denoted as $\{i\}_{1}^{\infty}$ such that
\begin{align}
J_{u_i}(r_i)> 2^{2d} \cdot J_{u_i}(\frac{r_i}{2})  \label{5.2.2}
\end{align}
where $r_i\in [r_0R_i, R_i]$.

Note $R_i^{-1}r_i\in [r_0, 1]$, without loss of generality, we can assume that there exists a subsequence of $\{i\}$, for simplicity also denoted as $\{i\}$ such that
\begin{align}
\lim_{i\rightarrow \infty}R_i^{-1}r_i= c_0\in [r_0, 1] \label{5.2.3}
\end{align}
where $c_0$ is some constant.

Taking the limit in (\ref{5.2.2}), from (\ref{5.2.3}) and Lemma \ref{lemma 5.1 construction of u_i} we get
\begin{align}
\frac{1}{\nu_{\infty}\big(B_{\infty}(c_0)\big)} \int_{B_{\infty}(c_0)} |u_{\infty}|^2 d\nu_{\infty}\geq \frac{2^{2d}}{\nu_{\infty}\big(B_{\infty}(\frac{c_0}{2})\big)} \int_{B_{\infty}\big(\frac{c_0}{2}\big)} |u_{\infty}|^2 d\nu_{\infty} \label{5.2.4}
\end{align}
From $u_{\infty}= r^{\alpha_1} \varphi_1(x)$ and $d> \alpha_1$, (\ref{5.2.4}) implies $\int_{X} |\varphi_1(x)|^2 dx= 0$, which is contradiction. 
}
\qed

Note $d\notin \mathscr{D}(M)$, from Theorem \ref{thm 3.6 induction of J} and induction method, there exists $k_0= k_0(d)$ such that for $r\in \big[\frac{k_0}{2}, R_i\big]$, (\ref{5.2.1}) holds, where we choose $i$ big enough such that $i\geq i_0(d- \alpha_1, r_0)$, $R_i> k_0$ and $u_i\nequiv 0$.

Now we define
\begin{align}
\check{u}_i(x)= \frac{u_i(x)}{\sqrt{J_{u_i}(\frac{k_0}{2})}} \label{def renormalized u_i}
\end{align}
then 
\begin{align}
J_{\check{u}_i}(\frac{k_0}{2})= 1\ , \quad\quad\quad \check{u}_i(p)= 0 \label{equation of renormalized u_i}
\end{align}

Note the scaling invariant property of (\ref{5.2.1}), hence there exists $i_0> 0$, if $i\geq i_0$, for $r\in [\frac{k_0}{2}, R_i]$, $J_{\check{u}_i}(r)\leq 2^{2d} \cdot J_{\check{u}_i}(\frac{r}{2})$, and we get 
\begin{align}
J_{\check{u}_i}(r)\leq \Big(2^{2d}\Big)^{\ln_2\big(\frac{r}{k_0}\big)+ 1} J_{\check{u}_i}(k_0)= \Big(\frac{r}{k_0}\Big)^{2d}\cdot \big(2^{2d}\big) J_{\check{u}_i}(k_0) \label{5.1}
\end{align}

By Theorem $1.2$ of \cite{LS},
\begin{align}
|\check{u}_i(x)|^2\leq C(n, p) J_{\check{u}_i}\big(2\rho(x)\big) \label{5.2}
\end{align}
recall that $\rho(x)= d(x, p)$.

From (\ref{5.1}) and (\ref{5.2}), when $i\geq i_0$, for $\rho(x)\in [\frac{k_0}{4}, \frac{R_i}{2}]$,
\begin{align}
|\check{u}_i(x)|\leq C(n, p, d, k_0) \rho(x)^{d} \label{5.3}
\end{align}

Combining with the Cheng-Yau's gradient estimate in \cite{CY} and the Arzela-Ascoli theorem, after taking suitable subsequence, $\check{u}_i$ converges to a polynomial growth harmonic function $u(x)$ on $M^n$. From (\ref{equation of renormalized u_i}), we know that $u(x)$ is not constant. The conclusion is proved.
}
\qed

{\it \textbf{Proof of Theorem \ref{thm existence of poly growth harmonic function}}}:~
{When the tangent cone at infinity of $M^n$ with renormalized limit measure is the unique metric cone with conic measure, denoted as $(C(X), \nu)$, then $\mathscr{D}(M)$ is a countable set by the fact that the spectrum of $(X, \nu_{-1})$ is a discrete set. Hence we can find $d\notin \mathscr{D}(M)$ and $d> \inf \{\alpha|\ \alpha\in \mathscr{D}(M), \alpha\neq 0\}$, from Theorem \ref{thm general existence of poly growth harmonic function}, the conclusion is proved.
}
\qed

\begin{lemma}\label{lem one criterion}
{Suppose $(M^n, g)$ has nonnegative sectional curvature, and for some fixed constants $\kappa> 1$, $a_0> 0$, 
\begin{align}
\lim_{s\rightarrow \infty}\frac{V(B(x, s))}{s^{\kappa}}= a_0 \label{oc 0}
\end{align}
where the convergence in (\ref{oc 0}) is uniform for all $x\in M^n$. Then the tangent cone at infinity of $M^n$ with renormalized limit measure is a unique metric cone $C(X)$ with unique conic measure $\nu$ of power $\kappa$, $\mathcal{H}^{\kappa- 1}(X)> 0$ and $\kappa\geq 2$ is an integer.
}
\end{lemma}

\begin{remark}\label{rem uni asy vol grow}
{If (\ref{oc 0}) holds uniformly for all $x\in M^n$, we will say $(M^n, g)$ has \textbf{uniform asymptotic polynomial volume growth of degree $\kappa$}.
}
\end{remark}

\pf
{Assume $x_i\rightarrow x$, $r_i\rightarrow \infty$, then $B_i(x_i, r)\rightarrow B_{\infty}(x, r)$, we have
\begin{align}
\nu\big(B_{\infty}(x, r)\big)&= \lim_{i\rightarrow \infty} \frac{\mu_i\Big(B_i(x_i, r)\Big)}{\mu_i\Big(B_i(p, 1)\Big)}= \lim_{i\rightarrow \infty}\frac{\mu\Big(B(x_i, r_ir)\Big)}{\mu\Big(B(p, r_i)\Big)} \nonumber \\
&=r^{\kappa} \lim_{i\rightarrow \infty}  \frac{\mu\Big(B(x_i, r_ir)\Big)}{a_0 (r_ir)^{\kappa}}= r^{\kappa} \label{oc 1}
\end{align}
where the last equation follows from the uniform convergence of (\ref{oc 0}).

From the definition of Hausdorff dimension, using (\ref{oc 1}), we obtain that the Hausdorff dimension of $C(X)$ is $\kappa$ and $\mathcal{H}^{\kappa}(C(X))> 0$. Because $C(X)$ is a metric cone on metric space $X$, it is not hard to get that the Hausdorff dimension of $X$ is $(\kappa- 1)$ and $\mathcal{H}^{\kappa- 1}(X)> 0$.

By Theorem $5.5$ of \cite{CC3} and the definitions of Ahlfors $l$-regular and $\nu$-rectifiable (Definition s $5.1$ and $5.3$ in \cite{CC3}), where $l$ is some non-negative number, using (\ref{oc 1}), we obtain that $\nu$ is Ahlfors $\kappa$-regular at all $x\in C(X)$, and $\kappa$ must be a non-negative integer. By assumption $\kappa> 1$, we proved that $\kappa\geq 2$ is an integer. 

From the Definition \ref{def conic measure} and (\ref{oc 1}), it is straightforward to verify that $\nu$ is a conic measure of power $\kappa$.


}
\qed

From the above Lemma and Theorem \ref{thm existence of poly growth harmonic function}, we have the following corollary.
\begin{cor}\label{cor section curvature case}
{Suppose $(M^n, g)$ has nonnegative sectional curvature and uniform asymptotic polynomial volume growth of degree $\kappa$, and $\kappa> 1$, then (\ref{upper bound of ds}) and (\ref{lower bound of dim}) hold.
}
\end{cor}

\begin{example}[Ding's example]\label{example Ding's}
{On $\mathbb{R}^{n}$, we define the warped product metric $g= dr^2+ f^2(r)d\mathbb{S}^{n- 1}$, where $\mathbb{S}^{n -1}$ is the classical $(n- 1)$-dimensional unit sphere, $f(r)$ is defined by modifying the famous symmetric mollifier $e^{-\frac{1}{1- r^2}}$ as the following :
\begin{equation}\label{Ding's f(r)}
f(r)= \left\{
\begin{array}{rl}
&a- b\exp\Big\{-\frac{1}{1- \big(r+ 3^{-\frac{1}{4}}\big)^2}\Big\}\ , \quad \quad 0\leq r< 1- 3^{-\frac{1}{4}} \\
&a \ , \quad \quad \quad \quad \quad \quad \quad \quad \quad r\geq 1- 3^{-\frac{1}{4}} 
\end{array} \right.
\end{equation}
where $b= a\cdot \exp\Big\{\frac{1}{1- 3^{-\frac{1}{2}}}\Big\}$, $a= \frac{\big(1- 3^{-\frac{1}{2}}\big)^2}{2\cdot 3^{-\frac{1}{4}}}$. The number $3^{-\frac{1}{4}}$ is chosen in (\ref{Ding's f(r)}), because it is the inflection point of the symmetric mollifier $e^{-\frac{1}{1- r^2}}$. It is straightforward to check that $Rc(g)\geq 0$ by the above definition of $f(r)$ and the metric $g$ is smooth. 

And it is obvious that $(\mathbb{R}^n, g)$ has linear volume growth and will not split isometrically. Hence by Theorem \ref{thm Sormani}, there does not exist any nonconstant harmonic function of polynomial growth on $(\mathbb{R}^n, g)$.
}
\end{example}

\begin{example}[Counterexample of Ni's Conjecture]\label{example Ni's}
{Let us start from the generalized Hopf fibration of $\mathbb{S}^7$ as the following:
\begin{align}
\mathbb{S}^3\longrightarrow \mathbb{S}^7\stackrel{\pi}{\longrightarrow} \mathbb{S}^4\ , \quad g^{\mathbb{S}^7}= k_1+ k_2 \nonumber
\end{align}
where $\mathbb{S}^3$, $\mathbb{S}^7$, $\mathbb{S}^4$ carry the metrics $g^{\mathbb{S}^3}$, $g^{\mathbb{S}^7}$, $\frac{1}{4}g^{\mathbb{S}^4}$; $\pi$ is a Riemannian submersion with totally geodesic fibers and $k_1= g^{\mathbb{S}^3}$, $k_2= \pi^{*}\big(\frac{1}{4}g^{\mathbb{S}^4}\big)$; $g^{\mathbb{S}^n}$ denotes the canonical metric of curvature $\equiv 1$ on $\mathbb{S}^n$.

Then for metric $g= dr^2+ f^2(r)k_1+ h^2(r)k_2$ on $M^8$, which is diffeomorphic to $\mathbb{R}^8$, from $(8.13)$ in \cite{CC1} and Section $2$ in \cite{BKN}, we have 
\begin{align}
Rm(X_1, X_2, X_1, X_2)&= \frac{1}{f^2(r)}- \Big(\frac{f'(r)}{f(r)}\Big)^2 \label{sc f}\\
Rm(X, Y, X, Y)&= \frac{f^2}{h^4}- \frac{f'\cdot h'}{f\cdot h}\label{sc f and h}\\
Rm(Y_1, Y_2, Y_1, Y_2)&= \frac{4}{h^2}- \frac{3f^2}{h^4}- \Big(\frac{h'}{h}\Big)^2 \label{sc h} \\
Rm(\partial r, X, \partial r, X)&= -\frac{f''}{f} \ , \quad \quad Rm(\partial r, Y, \partial r, Y)= -\frac{h''}{h} \label{sc r and f h}
\end{align}
where $X_1, X_2, X \in T\mathbb{S}^3$ and $Y_1, Y_2, Y\in T\mathbb{S}^4$.

In the following $a$, $\delta$, $\{c_i\}_{i= 0}^{6}$ are positive constants to be determined later, set 
\begin{equation}\label{Ni's f(r)}
f(r)= \left\{
\begin{array}{rl}
&\frac{1}{a}\sin(ar) \ , \quad \quad \quad \quad \quad \quad \quad \quad \quad 0\leq r\leq \delta \\
&\mathrm{f}(r- \delta) \ , \quad \quad \quad \quad \quad \quad \quad \quad \quad r> \delta
\end{array} \right.
\end{equation}
\begin{equation}\label{Ni's h(r)}
h(r)= \left\{
\begin{array}{rl}
&\frac{1}{a}\sin(ar) \ , \quad \quad \quad \quad \quad \quad \quad \quad \quad 0\leq r\leq \delta \\
&\mathrm{h}(r- \delta) \ , \quad \quad \quad \quad \quad \quad \quad \quad \quad r> \delta
\end{array} \right.
\end{equation}
where $\mathrm{f}$ and $\mathrm{h}$ are defined as the following:
\begin{align}
\mathrm{f}(x)\vcentcolon = c_1- c_2e^{-c_3x}\ , \quad \quad \mathrm{h}(x)\vcentcolon = c_0+ c_4x- c_5e^{-c_6x}\ , \quad \quad x\geq 0 \nonumber
\end{align}

It is easy to see that the metric $g= dr^2+ f^2(r)k_1+ h^2(r)k_2$ has positive sectional curvature when $0\leq r\leq \delta$, in the following we will try to find suitable constants $a$, $\delta$, $\{c_i\}_{i= 0}^{6}$ such that the $C^2$-metric $g$ has positive sectional curvature when $r> \delta$.

If $a$, $\delta$, $\{c_i\}_{i= 0}^{6}$ are positive constants satisfying the following eight assumptions:
 \begin{align}
a&= \Big(\frac{c_2c_3^2}{c_1- c_2}\Big)^{\frac{1}{2}}\ , \label{assp 1} \\
c_1&= c_2+ \frac{1- c_2^2c_3^2}{c_2c_3^2}\ , \label{assp 2} \\
\delta&= a^{-1}\sin^{-1}\big[a(c_1- c_2)\big]\ ,  \label{assp 3} \\
c_0- c_5&= c_1- c_2 \ , \label{assp 4}  \\ 
c_5c_6^2&= c_2c_3^2 \ , \label{assp 5} \\
c_4+ c_5c_6&= c_2c_3\ ,  \label{assp 6} \\
c_6&> c_3 \ ,        \label{assp 7} \\
c_0&\geq 3c_2+ c_5 \label{assp 8}
\end{align}
Note that there are many choices of $a$, $\delta$, $\{c_i\}_{i= 0}^{6}$ satisfying the above eight assumptions, and the following is one choice satisfying all the above assumptions:
\begin{align}
a&= \frac{1}{2\sqrt{3}} \ , \quad \delta= \frac{2\pi}{\sqrt{3}}\ , \quad c_0= \frac{13}{4}\ , \quad c_1= 4\ , \quad c_2= 1\ , \nonumber \\
c_3&= \frac{1}{2}\ , \quad c_4= \frac{1}{4}\ , \quad c_5= \frac{1}{4} \ , \quad c_6= 1 \nonumber 
\end{align}

Define $s(r)\vcentcolon = \frac{1}{a}\sin (ar)$, by (\ref{assp 1}), (\ref{assp 2}) and  (\ref{assp 3}), we have
\begin{align}
s(\delta)= \mathrm{f}(0)\ , \quad s'(\delta)= \mathrm{f}'(0)\ , \quad s''(\delta)= \mathrm{f}''(0) \label{f is c2}
\end{align}
which implies that $f(r)$ is a $C^2$ function on $[0, \infty)$. And by (\ref{assp 4}), (\ref{assp 5}) and (\ref{assp 6}), 
\begin{align}
\mathrm{f}(0)= \mathrm{h}(0)\ , \quad \mathrm{f}'(0)= \mathrm{h}'(0)\ , \quad \mathrm{f}''(0) = \mathrm{h}''(0) \label{h is c2}
\end{align}
And (\ref{h is c2}) combining with (\ref{f is c2}) yields that $h$ is a also a $C^2$ function on $[0, \infty)$.

From the definition of $f(r)$, $h(r)$ and the formula (\ref{sc r and f h}), it is easy to get 
\begin{align}
Rm(\partial r, X, \partial r, X)> 0\ , \quad \quad \quad Rm(\partial r, Y, \partial r, Y)> 0 \label{sc r and f h is posi}
\end{align}

Now we consider $(\mathrm{h}(x)- \mathrm{f}(x))''$, using (\ref{assp 7}), 
\begin{align}
(\mathrm{h}(x)- \mathrm{f}(x))''= c_2c_3^2e^{-c_6x}\big(e^{(c_6- c_3)x}- 1\big)\geq 0 \ , \quad \quad \forall x\geq 0 \label{ni's exap 1}
\end{align}
On the other hand, from (\ref{h is c2}), $(\mathrm{h}- \mathrm{f})'(0)= 0$. Then by (\ref{ni's exap 1})
\begin{align}
(\mathrm{h}- \mathrm{f})'(x)\geq (\mathrm{h}- \mathrm{f})'(0)= 0 \ , \quad \quad \forall x\geq 0 \label{ni's exap 2}
\end{align} 
Again, by (\ref{h is c2}) and (\ref{ni's exap 2}),
\begin{align}
(\mathrm{h}- \mathrm{f})(x)\geq (\mathrm{h}- \mathrm{f})(0)= 0 \ , \quad \quad \forall x\geq 0 \label{ni's exap 3}
\end{align} 

From (\ref{assp 2}), (\ref{assp 4}) and (\ref{assp 8}), we get
\begin{align}
\frac{1- c_2^2c_3^2}{c_2c_3^2}= c_1- c_2= c_0- c_5\geq 3c_2 \nonumber
\end{align}
simplifying it yields
\begin{align}
c_2c_3\leq \frac{1}{2} \label{ni's exap 4}
\end{align}
Then $\mathrm{h}'(0)= \mathrm{f}'(0)= c_2c_3\leq \frac{1}{2}$, by $\mathrm{h}''(x)< 0$, 
\begin{align}
\mathrm{h}'(x)< \mathrm{h}'(0)\leq \frac{1}{2} \label{ni's exap 5}
\end{align}

From (\ref{ni's exap 2}), (\ref{ni's exap 5}) and (\ref{sc f}), we obtain that when $r> \delta$,
\begin{align}
Rm(X_1, X_2, X_1, X_2)> 0 \label{sc f is posi}
\end{align}

From (\ref{ni's exap 3}), (\ref{ni's exap 5}) and (\ref{sc h}), when $r> \delta$,
\begin{align}
Rm(Y_1, Y_2, Y_1, Y_2)= \frac{\mathrm{h}^2\big(4- (\mathrm{h}')^2\big)- 3\mathrm{f}^2}{\mathrm{h}^4}> \frac{3(\mathrm{h}^2- \mathrm{f}^2)}{\mathrm{h}^4}\geq 0 \label{sc h is posi}
\end{align}

Now consider $\varphi(x)\vcentcolon= \mathrm{f}^3(x)- \mathrm{h}^3(x)\mathrm{f}'(x)\mathrm{h}'(x)$, note 
\begin{align}
\varphi(0)= \mathrm{h}^3(0)(1- \mathrm{f}'(0)^2)> 0 \label{ni's exap 6}
\end{align}
On the other hand, using $\mathrm{h}''< 0$, 
\begin{align}
\varphi'(x)&= 3\mathrm{f}' \mathrm{f}^2- 3\mathrm{h}^2(\mathrm{h}')^2\mathrm{f}'- \mathrm{h}^3 \mathrm{f}'' \mathrm{h}'- \mathrm{h}^3\mathrm{f}' \mathrm{h}'' \nonumber \\
&> \mathrm{h}'\mathrm{h}^2\big(-\mathrm{f}''\mathrm{h}-3\mathrm{f}'\mathrm{h}'\big) \label{ni's 7}
\end{align}
and 
\begin{align}
-\mathrm{f}''\mathrm{h}- 3\mathrm{f}'\mathrm{h}'&= \big(c_2c_3 e^{-c_3x}\big) \Big[c_0c_3+ c_3c_4x- 3c_4- (c_3c_5+ 3c_5c_6)e^{-c_6x}\Big] \nonumber \\
&\geq \big(c_2c_3 e^{-c_3x}\big) \Big[c_0c_3- 3c_4- (c_3c_5+ 3c_5c_6)\Big]  \nonumber \\
&= \big(c_2c_3 e^{-c_3x}\big) \Big[c_0c_3- 3c_2c_3- c_5c_3\Big]\geq 0 \label{ni's 8} 
\end{align}
in the last equation above we used (\ref{assp 6}), and in the last inequality we used (\ref{assp 8}). Combining (\ref{ni's 7}) with (\ref{ni's 8}), we obtain
\begin{align}
\varphi'(x)> 0 \label{ni's 9}
\end{align}
From (\ref{ni's exap 6}) and (\ref{ni's 9}), 
\begin{align}
\varphi(x)> 0 \ , \quad \quad \quad \quad \forall x\geq 0 \label{ni's 10}
\end{align}

By (\ref{ni's 10}) and (\ref{sc f and h}), when $r> \delta$, $x= r- \delta> 0$, 
\begin{align}
Rm(X, Y, X, Y)= \frac{\varphi(x)}{\mathrm{f}(x)\cdot \mathrm{h}^4(x)}> 0 \label{sc f and h is posi}
\end{align}

From (\ref{sc r and f h is posi}), (\ref{sc f is posi}), (\ref{sc h is posi}) and (\ref{sc f and h is posi}), the metric $g= dr^2+ f^2(r)k_1+ h^2(r)k_2$ on $M^8$ has positive sectional curvature, where $f$, $g$ are defined in (\ref{Ni's f(r)}) and (\ref{Ni's h(r)}). It is not hard to see that this metric also has the uniform asymptotic polynomial volume growth of degree $5$ as in (\ref{oc 0}) and Remark \ref{rem uni asy vol grow}. Then by Corollary \ref{cor section curvature case}, there exists nonconstant harmonic function of polynomial growth on $(M^8, g)$, but $(M^8, g)$ does not have maximal volume growth. This disproves the necessary part of Conjecture \ref{conj 5.4}.
}
\end{example}

\section{Uniform bound of frequency function}\label{SECTION 5}

Much of argument in this section followed the detailed analysis about frequency function in \cite{CM-JDG-97} (especially Proposition $3.3$, Proposition $3.36$, Proposition $4.11$ and Lemma $7.1$ there). We are providing details here again to make our argument self-contained and concrete enough for our purpose, some more general argument can be found in \cite{CM-JDG-97}.

In this section, $I(r)$, $D(r)$, $E(r)$, $F(r)$, $\mathscr{F}(r)$ and $\mathscr{W}(r)$ are defined as in Section \ref{SECTION 2} with respect to some nonconstant function $u\in \mathscr{H}_d(M)$. Further assume that $u(p)= 0$, where $p\in M^n$.

The following Lemma is a weak version of a uniform Harnack inequality for harmonic function with polynomial growth.

\begin{lemma}\label{lemma 4.1 weak Harnack inequality}
{For nonconstant $u\in \mathscr{H}_d(M)$, there exists positive increasing sequence $\{r_i\}_{i= 1}^{\infty}$ such that $\lim_{i\rightarrow \infty}r_i= \infty$ and for any $i$
\begin{align}
D(2^{4n+ 1}r_i)\leq 2^{10nd}D(r_i) \nonumber
\end{align}
}
\end{lemma}

\pf
{By contradiction. If the conclusion does not hold, there exists $R_0\geq 1$, such that
\begin{align}
D(2^{4n+ 1}r)> 2^{10nd}D(r) \quad \quad \quad when \ r\geq R_0 \nonumber
\end{align}
by induction we get that for any $j= 1, 2, 3, \cdots$
\begin{align}
D(2^{(4n+ 1)j} R_0)> 2^{10ndj}D(R_0) \label{4.1.1}
\end{align}
 
On the other side, by $u\in \mathscr{H}_d(M)$ and Corollary $3.2$ in Chapter $1$ of \cite{SY},
\begin{align}
\sup_{B(r)}|\nabla u|\leq \frac{C(n)}{r}\sup_{B(2r)} |u|\leq C(n, d)r^{d- 1}\ , \quad \quad \forall r\geq 1 \label{bound of gradient}
\end{align}

From (\ref{estimate of b}) and $V(r)\leq V_0^n(1) r^n$, there exists $C_1> 0$ such that for any $r> 0$,  
\begin{align}
D(r)= r^{2- n}\int_{b(x)\leq r} |\nabla u(x)|^2 dx\leq C_1 (r^{2d}+ 1) \label{4.1.2}
\end{align}

By (\ref{4.1.1}) amd (\ref{4.1.2}),
\begin{align}
C_1\Big[\big(2^{(4n+ 1)j}R_0\big)^{2d}+ 1\Big]> 2^{10ndj}D(R_0) \nonumber
\end{align}
which implies
\[2^{(2- 2n)jd}+ 2^{-10ndj}> \frac{D(R_0)}{C_1 R_0^{2d}}\]
let $j\rightarrow \infty$ in the above, we get $0\geq \frac{D(R_0)}{C_1 R_0^{2d}}$. However, $D(R_0)> 0$ because $u$ is nonconstant, which is the contradiction. 
}
\qed

\begin{lemma}[Equivalence of $E$ and $D$]\label{lemma 4.2 the ratio between D and E}
{For $\epsilon> 0$, there exists $\mathfrak{R}_1= \mathfrak{R}_1(\epsilon, p, n, V_M, \gamma)> 0$ such that for $r\geq \mathfrak{R}_1$, if 
\begin{align}
D(2^{4n+ 1}r)\leq \gamma D(r) \label{4.2.1}
\end{align}
then for $s\in [r, 2^{4n}r]$, 
\begin{align}
\ln \Big(\frac{D(s)}{E(s)}\Big)\leq \epsilon \label{4.2.2}
\end{align}
}
\end{lemma}

\pf
{From Theorem \ref{thm collecting facts}, for given $\delta> 0$, there exists $\mathfrak{R}= \mathfrak{R}(p, \delta)> 0$ such that for $\rho(x)=r\geq \mathfrak{R}$,
\begin{align}
\Big|\ln \frac{b(x)}{\rho(x)}\Big|\leq \delta\ , \quad \quad \int_{b(x)\leq r} \Big||\nabla b|^2- 1\Big|^2 dx\leq \delta^2 \mathrm{Vol}(b(x)\leq r)  \label{4.2.3}
\end{align}
where $\delta> 0$ is to be determined later. From Cauchy-Schwarz inequality,
\begin{align}
\int_{b\leq r} \Big||\nabla b|^2- 1\Big|\leq \delta \mathrm{Vol}(b\leq r) \label{4.2.4}
\end{align}

Then for $s\in [r, 2^{4n}r]$, 
\begin{align}
|D(s)- E(s)|&= s^{2- n}\int_{b\leq s} |\nabla u|^2 (1- |\nabla b|^2) \nonumber \\
&\leq s^{2- n} \sup_{b(x)\leq s} |\nabla u|^2(x)\cdot \delta \mathrm{Vol}(b\leq s) \nonumber \\
&\leq C(\delta, n) s^2 \sup_{b(x)\leq s} |\nabla u|^2(x) \label{4.2.5}
\end{align}
in the first inequality above we used (\ref{4.2.4}), and we have $\lim_{\delta\rightarrow 0} C(\delta, n)= 0$.

Without loss of generality, assume $e^{\delta}\leq \sqrt{\frac{4}{3}}$. From Theorem $1.2$ in \cite{LS} and (\ref{4.2.3}) above, 
\begin{align}
\sup_{b(x)\leq 2^{4n}r} |\nabla u|^2(x)&\leq \sup_{B(\sqrt{\frac{4}{3}}\cdot 2^{4n}r)} |\nabla u|^2\leq \frac{C(n)}{\mathrm{Vol}(B(\sqrt{3}\cdot 2^{4n}r))} \int_{B(\sqrt{3}\cdot 2^{4n}r))} |\nabla u|^2 \nonumber \\
&\leq \frac{C(n)}{V_M\cdot \Big(\sqrt{3}\cdot 2^{4n}r\Big)^n} \int_{b\leq 2^{4n+1}r} |\nabla u|^2 \nonumber \\
&\leq C(n, V_M)r^{-2} D(2^{4n+ 1}r)\leq C(n, V_M)\gamma D(r)r^{-2}  \label{4.2.6}
\end{align}
in the last inequality we used (\ref{4.2.1}).

Note that $s\in [r, 2^{4n}r]$, hence 
\begin{align}
D(r)\leq r^{2- n}\int_{b\leq s} |\nabla u|^2\leq \Big(\frac{s}{r}\Big)^{n- 2} D(s)\leq 2^{4n(n- 2)} D(s) \label{4.2.7}
\end{align}

From (\ref{4.2.6}) and (\ref{4.2.7}), 
\begin{align}
\sup_{b(x)\leq 2^{4n}r} |\nabla u|^2(x)\leq C(n, V_M)\gamma D(s)r^{-2}\leq C(n, V_M)\gamma D(s) s^{-2} \label{4.2.8}
\end{align}
in the last inequality we used $s\leq 2^{4n}r$.

By (\ref{4.2.5}) and (\ref{4.2.8}), we obtain
\begin{align}
|D(s)- E(s)|\leq C(\delta, n, V_M)\gamma D(s) \nonumber
\end{align}
where $\lim_{\delta\rightarrow 0} C(\delta, n, V_M)= 0$. Hence 
\begin{align}
\frac{E(s)}{D(s)}\geq 1- \gamma C(\delta, n, V_M) \label{4.2.9}
\end{align}

There exists $\delta= \delta(n, V_M, \epsilon, \gamma)$ such that (\ref{4.2.9}) implies $\ln \Big(\frac{D(s)}{E(s)}\Big)\leq \epsilon$. Combining all the above, there exists $\mathfrak{R}_1= \mathfrak{R}_1(p, \delta)= \mathfrak{R}_1(p, n, V_M, \epsilon, \gamma)$ satisfying our conclusion.
}
\qed

\begin{lemma}\label{lem 4.3 almost monotonicity of W}
{Given positive constants $\gamma$, $\epsilon$, there exists $\mathfrak{R}_2= \mathfrak{R}_2(\epsilon, p, n, V_M, \gamma)> 0$ such that if $r> \mathfrak{R}_2$ and 
\begin{align}
D(2^{4n+ 1}r)\leq \gamma D(r) \label{4.3.1}
\end{align}
then
\begin{align}
\int_{r}^{2^{4n}r} \min\{(\ln \mathscr{W})'(t), 0\} dt> -\epsilon \label{4.3.2}
\end{align}
}
\end{lemma}

\pf
{Using the first variation formula of energy in the Appendix of \cite{CM-JDG-97}, we have 
\begin{align}
\Big(\ln \mathscr{W}(s)\Big)'&= \frac{2}{s}+ \frac{2F(s)}{E(s)}- \frac{s^{1- n}\int_{b\leq s} Hess(b^2) (\nabla u, \nabla u)}{E(s)}- \frac{2D(s)}{sI(s)} \nonumber \\
&= \mathfrak{W}+ \mathfrak{J}+ \mathfrak{K} \label{4.3.3}
\end{align}
where 
\begin{align}
\mathfrak{W}&= \frac{2s^{2- n}\int_{b= s} \big|\frac{\partial u}{\partial n}\big|^2 |\nabla b|^{-1}}{D(s)}- \frac{2D(s)}{sI(s)} \ , \quad \quad \mathfrak{J}= \frac{2}{s}- \frac{s^{1- n}\int_{b\leq s} Hess(b^2) (\nabla u, \nabla u)}{E(s)}\ , \nonumber\\
\mathfrak{K}&= \frac{2F(s)}{sE(s)}- \frac{2s^{2- n}\int_{b= s} \big|\frac{\partial u}{\partial n}\big|^2 |\nabla b|^{-1}}{D(s)} \nonumber
\end{align}

From Cauchy-Schwarz inequality, it is easy to get $\mathfrak{W}\geq 0$. Hence we only need to bound the integrals of $\mathfrak{J}$ and $\mathfrak{K}$.

From Theorem \ref{thm collecting facts}, given any $\delta> 0$, there exists $\mathfrak{R}= \mathfrak{R}(p, \delta)> 0$ such that for $r= \rho(x)\geq \mathfrak{R}$, we have
\begin{align}
\Big|\ln \frac{b(x)}{\rho(x)}\Big|&\leq \delta \ , \label{4.3.4} \\
\int_{b(y)\leq r} \Big||\nabla b|^2- 1\Big|^2 dy &\leq \delta^2 \mathrm{Vol}(b\leq r) \ , \label{4.3.5} \\
\int_{b(y)\leq r} \Big|Hess(b^2)- 2g\Big|^2 dy &\leq \delta^2 \mathrm{Vol}(b\leq r) \label{4.3.6}
\end{align}
We choose $\delta> 0$ such that 
\begin{align}
e^{2\delta}< \frac{4}{3}  \label{4.3.7}
\end{align}
then (\ref{4.3.4}) implies that for $s\geq \mathfrak{R}$, 
\begin{align}
\{b\leq s\}\subset B\Big(\sqrt{\frac{4}{3}}s\Big)\ , \quad\quad \quad B(\sqrt{3}s)\subset \{b\leq 2s\} \nonumber
\end{align}

Now we estimate $\mathfrak{J}$,
\begin{align}
|\mathfrak{J}|&\leq \frac{2}{s}\Big[\frac{D(s)- E(s)}{E(s)}\Big]+ \frac{1}{s}\cdot \frac{s^2\sup_{b\leq s}|\nabla u|^2}{E(s)}\cdot \Big[s^{-n}\int_{b\leq s} \Big|Hess(b^2)- 2g\Big|\Big] \nonumber \\
&= \mathfrak{J}_1+ \mathfrak{J}_2 \nonumber
\end{align}

From Lemma \ref{lemma 4.2 the ratio between D and E}, there exists $\mathfrak{R}_1= \mathfrak{R}_1(\epsilon, p, n, V_M, \gamma)$ such that if $r\geq \max\big\{\mathfrak{R}, \mathfrak{R}_1\big\}$, for $s\in [r, 2^{4n}r]$, $\ln\frac{D(s)}{E(s)}\leq \delta$, hence $\mathfrak{J}_1\leq \frac{2}{s}\Big[e^{\delta}- 1\Big]$. 

From (\ref{4.2.8}), (\ref{4.3.6}), the Cauchy-Schwarz inequality and the Bishop-Gromov volume comparison theorem, we get
\begin{align}
\mathfrak{J}_2\leq \frac{1}{s}C(n, V_M)\gamma \frac{D(s)}{E(s)}\cdot \delta V_0^n(1)e^{n\delta}\leq \frac{1}{s}C(n, V_M, \gamma)\delta e^{[(n+ 1)\delta]}\leq \frac{C(n, V_M, \gamma)}{s}\delta \nonumber
\end{align}
in the last inequality we used (\ref{4.3.7}).

Hence $|\mathfrak{J}|\leq \frac{2}{s}\Big[e^{\delta}- 1\Big]+ \frac{C(n, V_M, \gamma)}{s}\delta$, taking integral on $[r, 2^{4n}r]$, where $r\geq \max\big\{\mathfrak{R}, \mathfrak{R}_1\big\}$, 
\begin{align}
\int_{r}^{2^{4n}r} |\mathfrak{J}|\leq \Big\{2\big[e^{\delta}- 1\big]+ C(n, V_M, \gamma)\delta \Big\}\ln (2^{4n}) \label{4.3.8}
\end{align}

Next we estimate the integral of $\mathfrak{K}$,
\begin{align}
\mathfrak{K}&= \frac{2}{s}\Big[\frac{D(s)- E(s)}{E(s)}\Big]\cdot \frac{F(s)}{D(s)}+ \frac{2s^{2- n}\int_{b= s} \big|\frac{\partial u}{\partial n}\big|^2\big(|\nabla b|- |\nabla b|^{-1}\big)}{D(s)} \nonumber \\
&= \mathfrak{K}_1+ \mathfrak{K}_2 \nonumber 
\end{align}

From \ref{estimate of b}, it is easy to see $D(s)\geq E(s)$, hence $\mathfrak{K}\geq \mathfrak{K}_2$. Now,
\begin{align}
|\mathfrak{K}_2|&\leq 2s^{2- n} \frac{\sup_{b= s} |\nabla u|^2}{D(s)} \int_{b(x)= s} \Big| |\nabla b|- |\nabla b|^{-1}\Big|dx \nonumber \\
&\leq 2s^{2- n} \cdot \frac{\Big(\sup_{b\leq 2^{4n}r} |\nabla u|^2\Big)}{D(s)} \Big| |\nabla b|- |\nabla b|^{-1}\Big|dx \nonumber \\
&\leq C(n, V_M, \gamma)s^{-n}\int_{b= s} \Big| |\nabla b|- |\nabla b|^{-1}\Big| \nonumber
\end{align}
we used (\ref{4.2.8}) in the last inequality.

Hence
\begin{align}
\int_{r}^{2^{4n}r} \mathfrak{K} ds&\geq -C(n, V_M, \gamma)\int_{b\leq 2^{4n}r} r^{-n} \big||\nabla b|^2- 1\big| \nonumber \\
&\geq -C(n, V_M, \gamma) \frac{\mathrm{Vol}(b\leq 2^{4n}r)}{r^n}\delta e^{n\delta}\geq -C(n, V_M, \gamma)\delta  \label{4.3.10}
\end{align}
in the first inequality above we used the co-area formula, and (\ref{4.3.5}) was used in the second inequality. From (\ref{4.3.8}) and (\ref{4.3.10}), 
\[\lim_{\delta\rightarrow \infty} \int_{r}^{2^{4n}r} \big(\mathfrak{J}+ \mathfrak{K}\big)\geq \lim_{\delta\rightarrow \infty} -C(n, V_M, \gamma)\big[e^{\delta}- 1+ \delta \big]= 0 \]

On the other side, we have
\begin{align}
\int_{r}^{2^{4n}r} \min\{\big(\ln \mathscr{W}\big)'(t), 0\} dt\geq \int_{r}^{2^{4n}r} \big(\mathfrak{J}+ \mathfrak{K}\big) \nonumber
\end{align}

Hence there exists $\delta_0= \delta_0(\epsilon, n, V_M, \gamma)$ satisfying (\ref{4.3.7}), and if $\delta\leq \delta_0$,
\begin{align}
\int_{r}^{2^{4n}r} \min\{\big(\ln \mathscr{W}\big)'(t), 0\} dt\geq -\epsilon  \nonumber
\end{align}

Choose $\mathfrak{R}_2= \max\{\mathfrak{R}(p, \delta_0), \mathfrak{R}_1(\epsilon, p, n, V_M, \gamma)\}$, the conclusion is proved.
}
\qed

\begin{lemma}\label{lem 4.4 bounding growth of D by growth of I}
{For $p\in M$, there exists $\mathfrak{R}_3= \mathfrak{R}_3(p)> 0$ such that if $r> \mathfrak{R}_3$ and 
\begin{align}
I(2^{4n+ 2}r)\leq \gamma I(\frac{r}{2}) \label{4.4.1}
\end{align}
then
\begin{align}
D(2^{4n+ 1}r)\leq C_1(n)\gamma D(r) \label{4.4.2}
\end{align}
where $C_1(n)$ is the constant depending only on $n$.
}
\end{lemma}

\pf
{As in the proof of Lemma \ref{lemma 4.2 the ratio between D and E}, we can choose $\mathfrak{R}_3= \mathfrak{R}_3(p)> 0$ such that for $\rho(x)=r > \frac{\mathfrak{R}_3}{2}$, $|\ln \frac{b(x)}{\rho(x)}|\leq \frac{1}{2}\ln \frac{4}{3}$. Similar as (\ref{4.2.6}), for $r> \mathfrak{R}_3$,
\begin{align}
\sup_{b(x)\leq \frac{r}{2}} |\nabla u|^2(x)\leq \frac{C(n)}{V_M}r^{-2} D(r) \label{4.4.3}
\end{align}

Integrating (\ref{4.4.3}) along geodesics starting at $p$ and using $u(p)= 0$, we obtain 
\begin{align}
\sup_{b(x)\leq \frac{r}{2}} |u|^2(x)\leq \frac{4}{3}\frac{C(n)}{V_M} D(r) \nonumber
\end{align}

Hence
\begin{align}
I(\frac{r}{2})&= \big(\frac{r}{2}\big)^{1- n}\int_{b= \frac{r}{2}} u^2|\nabla b|\leq \frac{C(n)}{V_M} D(r)\cdot \big(\frac{r}{2}\big)^{1- n} \int_{b= \frac{r}{2}} |\nabla b|\nonumber \\
&= \frac{C(n)}{V_M} D(r) I_1(r) = C(n) D(r) \label{4.4.4}
\end{align}
in the last equation we used (\ref{I_1(r)}).

From (\ref{I'(r)}),
\begin{align}
\int_{2^{4n+ 1}r}^{2^{4n+ 2}r} \frac{2D(s)}{s} ds= I(2^{4n+ 2}r)- I(2^{4n+ 1}r)\leq I(2^{4n+ 2}r) \nonumber 
\end{align}
which implies
\begin{align}
2\int_{2^{4n+ 1}r}^{2^{4n+ 2}r} s^{n- 2}D(s) ds\leq \Big(2^{4n+ 2}r\Big)^{n- 1}I(2^{4n+ 2}r) \nonumber
\end{align}

Note $s^{n- 2}D(s)$ is nondecreasing in $s$ from the definition of $D(r)$, we get
\begin{align}
2\Big(2^{4n+ 1}r\Big)^{n- 1} D(2^{4n+ 1}r)\leq \Big(2^{4n+ 2}r\Big)^{n- 1}I(2^{4n+ 2}r) \nonumber
\end{align}

Combining (\ref{4.4.4}), simplifying the above inequality yields 
\begin{align}
D(2^{4n+ 1}r)\leq 2^{n- 2} I(2^{4n+ 2}r)\leq 2^{n -2}\gamma I(\frac{r}{2})= C_1(n)\gamma D(r) \nonumber
\end{align}
}
\qed

\begin{theorem}\label{thm 4.5 frequency of linear growth function}
Suppose that $(M^n, g)$ has nonnegative Ricci curvature and maximal volume growth. For $u\in \mathscr{H}_d(M)$, if $5d\notin \mathscr{D}(M)$, then the frequency of $u$ is bounded by $C(u, n, V_M, d)$.
\end{theorem}

\pf
{If $u$ has no zero point, then by Yau's Liouville theorem \cite{Yau} $u$ is constant, the conclusion is straightforward. Assume $u(p)= 0$ where $p\in M$ is some fixed point. We will firstly prove the following claim:
\begin{claim}\label{claim 4.6}
{There exists a constant $\mathfrak{R}= \mathfrak{R}(p, n, V_M, d)> 0$ such that if $r\geq \mathfrak{R}$, $\mathscr{F}(r)\leq C(n, d)$.
}
\end{claim}

By Lemma \ref{lemma 4.1 weak Harnack inequality}, there exist $r_i\rightarrow \infty$, such that
\begin{align}
D(2^{4n+ 1}r_i)\leq 2^{10nd}D(r_i) \nonumber
\end{align}

Choose $\delta> 0$, such that $\delta< \frac{1}{2}\ln \frac{4}{3}$, by Lemma \ref{lemma 4.2 the ratio between D and E} and Lemma \ref{lem 4.3 almost monotonicity of W}, there exists $\mathfrak{R}_4= \mathfrak{R}_4(\delta, p, n, V_M, 2^{10nd})> 0$ such that if $r_i> \mathfrak{R}_4$, then
\begin{align}
&\ln \Big(\frac{D(s)}{E(s)}\Big)\leq \delta \ , \quad\quad\quad \forall s \in [r_i, 2^{4n}r_i] \label{4.5.1} \\
&\int_{r_i}^{2^{4n}r_i} \min\{(\ln \mathscr{W})'(t), 0\}dt > -\delta \label{4.5.2} 
\end{align}

From (\ref{I'(r)}), we get 
\begin{align}
\big(\ln I(r)\big)'= \frac{2\mathscr{F}(r)}{r} \label{the derivative of ln(I)}
\end{align}

Hence 
\begin{align}
&\int_{r_i}^{2^{4n}r_i} \frac{2\mathscr{F}(s)}{s} ds= \Big(\ln I(s)\Big)\Big|_{r_i}^{2^{4n}r_i}= \ln\Big(\frac{D(2^{4n}r_i)}{D(r_i)}\Big)- \ln\Big(\frac{\mathscr{F}(2^{4n}r_i)}{\mathscr{F}(r_i)}\Big) \nonumber \\
&\quad  \leq \ln \Big(\frac{D(2^{4n}r_i)}{D(2^{4n+ 1}r_i)}\cdot 2^{10nd}\Big)+ \ln \Big(\frac{E(2^{4n}r_i)}{D(2^{4n}r_i)}\Big)+ \ln \Big(\frac{D(r_i)}{E(r_i)}\Big)- \ln \Big(\frac{\mathscr{W}(2^{4n}r_i)}{\mathscr{W}(r_i)}\Big) \nonumber \\
&\quad \leq (10nd+ n- 2)\ln 2+ 3\delta\leq 12dn\cdot \ln 2 \label{4.5.3}
\end{align}
in the second inequality from the end, we used (\ref{4.5.1}) and (\ref{4.5.2}).

From (\ref{4.5.3}), there exists $s_i\in [2^{2n}r_i, 2^{4n}r_i]$ such that 
\begin{align}
\mathscr{F}(s_i)\leq 3d \label{4.5.4}
\end{align}

By (\ref{4.5.1}), (\ref{4.5.2}) and (\ref{4.5.4}), for $r_i> \mathfrak{R}_4$, $s\in [r_i, 2^{2n}r_i]$,
\begin{align}
\mathscr{W}(s)\leq 3d e^{2\delta}\leq 4d \nonumber
\end{align}
Combining with (\ref{4.5.1}), we get
\begin{align}
\mathscr{F}(s)\leq \mathscr{W}(s)e^{\delta}< 5d \label{bound of F}
\end{align}
where $r_i> \mathfrak{R}_4$ and $s\in [r_i, 2^{2n}r_i]$.

Using (\ref{the derivative of ln(I)}), for $r\in [2r_i, 2^{2n}r_i]$, we get 
\[I(r)\leq 2^{10d} I(\frac{r}{2})\]

By $5d\notin \mathscr{D}(M)$, using Theorem \ref{thm 3.5 induction of I} and induction method, there exists $k_0= k_0(d)$ such that if $r_i> k_0$ then 
\begin{align}
I(r)\leq 2^{10d} I(\frac{r}{2}) \ , \quad \quad \quad r\in [\frac{k_0}{2}, 2^{2n}r_i] \nonumber
\end{align}
which implies $I(2^{4n+ 2}r)\leq \big(2^{10d}\big)^{4n+ 3}I(\frac{r}{2})$ for $r\in [\frac{k_0}{2}, 2^{-2-2n}r_i]$.

Let $\mathfrak{R}_5= \max\{\frac{k_0}{2}, \mathfrak{R}_3\}$, where $\mathfrak{R}_3$ is from Lemma \ref{lem 4.4 bounding growth of D by growth of I}. Note $\mathfrak{R}_5= \mathfrak{R}_5(p, d)$. By Lemma \ref{lem 4.4 bounding growth of D by growth of I}, 
\begin{align}
D(2^{4n+1}r)\leq C_1(n)2^{60nd} D(r) \ , \quad\quad \quad r\in [\mathfrak{R}_5, 2^{-2-2n}r_i] \label{4.5.5}
\end{align}

From (\ref{4.5.5}), similar to the above argument to get (\ref{bound of F}), we get that there exists $\mathfrak{R}= \max\Big\{\mathfrak{R}_5, \mathfrak{R}_4\big(\delta, p, n, V_M, C_1(n)2^{60nd}\big)\Big\}$ such that for $r\in [\mathfrak{R}, 2^{-2-2n}r_i]$, 
\begin{align}
\mathscr{F}(s)\leq 31d+ C_1(n) \ , \quad\quad \quad s\in [\mathfrak{R}, 2^{-2}r_i] \nonumber
\end{align}

In the above inequality, let $i\rightarrow \infty$, then for $r\geq \mathfrak{R}(\delta, p, n ,V_M, d)$, $\mathscr{F}(r)\leq C(n, d)$. If we fix $\delta= \frac{1}{4}\ln\frac{4}{3}$, then Claim \ref{claim 4.6} is proved.

Because $\mathscr{F}(r)$ is continuous function of $r$, $\mathscr{F}(r)\leq C$ on $[0, \mathfrak{R}]$, where $C$ is some constant depending on $u$ and $\mathfrak{R}$. Combining the above results together, the conclusion of the theorem is proved.
}
\qed

Now we prove Theorem \ref{thm unique cone case} by using the above theorem.

{\it \textbf{Proof of Theorem \ref{thm unique cone case}}}:~
{If $d< 1$, from \cite{Cheng-sub} $u$ must be constant, then the conclusion follows trivially. Hence we assume $d\geq 1$ in the rest of the proof. When the tangent cone at infinity of $M^n$ is unique, denoted as $C(X)$, then $\mathscr{D}(M)$ is a countable set by the fact that the spectrum of $C(X)$ is discrete.  Because $d\geq 1$, we can find $d_0\in [d, d+ 1]$ such that $5d_0\notin \mathscr{D}(M)$, note $u\in \mathscr{H}_d(M)\subset \mathscr{H}_{d_0}(M)$. From Theorem \ref{thm 4.5 frequency of linear growth function}, $\mathscr{F}_u(r)\leq C(u, n, V_M, d_0)\leq C(u, n, V_M, d)$, the conclusion is proved.
}
\qed

\section*{Acknowledgments}
The author was partially supported by NSFC 11401336. We thank Jiaping Wang for his interest and continuous encouragement, Xian-Tao Huang, William P. Minicozzi II, Christina Sormani, Shing-Tung Yau for their comments, and Liqun Zhang for sending the offprint \cite{Zhang} to us. We are indebted to Bo Yang for his comments and pointing out the relation between Conjecture \ref{conj 5.4} and the result in \cite{Ding2} to us in 2012. We are grateful to Shouhei Honda for his detailed comments and enthusiastic suggestions on the paper. Last but not least, we particularly thank Gang Liu for carefully reading the earlier version of the paper and pointing out some gaps, and we benefit from several long conversations with him.


\end{document}